\documentclass[a4paper,12pt]{article}
\usepackage{amsmath,amssymb,theorem}
\usepackage{amscd}
\usepackage{graphicx}

\newcommand{\bma}{\begin{math}}
\newcommand{\ema}{\end{math}}
\newcommand{\beq}{\begin{eqnarray}}
\newcommand{\eeq}{\end{eqnarray}}
\newcommand{\be}{\begin{eqnarray*}}
\newcommand{\ee}{\end{eqnarray*}}

\theorembodyfont{\itshape}
\newtheorem{thm}{Theorem}
\newtheorem{prop}[thm]{Proposition}
\newtheorem{lm}[thm]{Lemma}

\theorembodyfont{\rmfamily}
\newtheorem{df}[thm]{Definition}
\newtheorem{pf}{Proof.}
\newtheorem{re}{Remark}

\newtheorem{N}{Notation}

\newcommand{\qed}{\hfill $\Box$ \par}

\begin{document}

\title{ Generalization  of Schl\"afli formula to the volume of a spherically
faced simplex}
\author{Kazuhiko AOMOTO and Yoshinori MACHIDA}

\date{}

\maketitle

\begin{abstract}

The purpose of this note is to present two identities (contiguity relation
and variation formula) concerning the volume of a spherically faced simplex
in the Euclidean space.
These identities are described in terms of Cayley-Menger determinants
and their differentials involved with hypersphere arrangements.
They are derived as a limit of fundamental identities for
hypergeometric integrals associated with hypersphere arrangements
obtained  by the authors in the preceding article.
The corrected version of the variational formula of Schl\"afli type for the volume of a spherically faced simplex
is presented.

\end{abstract}




\renewcommand{\thefootnote}{\fnsymbol{footnote}}



\footnotetext{\scriptsize
{\it Running Title}: Variation formula for the volume of a spherically faced simplex.}

\bigskip

\section{Introduction and preliminaries}



\bigskip

In this article, we give a new variation formula for the volume of a pseudo simplex 
with spherical faces in the Euclidean space (see Theorems I and II). 
To derive it, we apply the variation formula
obtained in [8], which is involved in  hypergeometric integrals associated 
with hypersphere arrangements. 
This procedure can be done by regularization of integrals (the method of generalized functions)
i.e., by taking the zero limit of exponents for hypergeometric integrals (see [11]).
A hypersphere arrangement in the $n$ dimensional  Euclidean space can be realized 
by the stereographic projection 
as the restriction to  the fundamental unit hypersphere  of a hyperplane arrangement 
in the $(n+1)$ dimensional Euclidean space.  
The theory of
hypergeometric integrals associated with hypersphere arrangements 
has been developed in this framework
in terms of twisted rational de Rham cohomology  (see [6], [8]). 
It is described in terms of  Cayley-Menger determinants.
In Theorem III, we give a correct version of some errors in the variation volume formula in [4],
which is  an extension of 
 L.Schl\"afli formula of a geodesic simplex in the unit hypersphere
 (see [1], [2], [3], [9], [12], [15], [16], [20], [22]. [23], 
and see also [9], [20], [21] related to the bellows conjecture).


\bigskip

1. Let $n+1$ real quadratic polynomials $f_j$ in $n$ variables
$x = (x_1,\ldots,x_n)$ in ${\bf R}^n$
\be
f_j(x)  = Q(x) + \sum_{\nu=1}^n 2\alpha_{j\nu}x_\nu + \alpha_{j0}\quad(1\le j\le n+1)
\ee
be given, where $Q(x)$ denotes the quadratic form
\be
Q(x) = \sum_{j=1}^n x_j^2,
\ee
and $\alpha_{j\nu}\in {\bf R},\, \alpha_{j0} \in {\bf R}$.
 
Let $S_j$ be the $n-1$ dimensional hypersphere defined by $f_j(x) = 0.$ 
Denote by $O_j$ the center of $S_j$, by $r_j$ the radius of $S_j$
and by $\rho_{jk}$ the distance of $O_j$ and $O_k.$

Then 
\be
&& r_j^2 = - \alpha_{j0} + \sum_{\nu=1}^n \alpha_{j\nu}^2,\\
&& \rho_{jk}^2 = \sum_{\nu = 1}^n (\alpha_{j\nu} - \alpha_{k\nu})^2.
\ee

Let ${\cal A} = \bigcup_{j=1}^{n+1} S_j$ be the arrrangement of hyperspheres 
consisting of $n+1$  hyperspheres $S_j.$

\medskip

2. Denote by $N$ the set of indices $\{1,2,\ldots,\,n+1\}$.

\begin{df}
Denote by $J, K$ the two non-empty subsets of indices
$J  = \{j_1,\ldots, j_p\}, K = \{k_1,\ldots,k_p\} \subset N$
of size $p.$
Cayley-Menger determinants associated with $\cal A$ are given by 
the following system of determinants (see [9], [12], [15], [25]):
\be
&&
B\!
\left(
\begin{array}{ccccc}
 0 &j_1   &\ldots& j_p   \\
 0 &k_1   &\ldots&k_p   \\
\end{array}
\right) : =
\left|
\begin{array}{cccccc}
  0&1   &1&\ldots&1   \\
  1& \rho_{j_1k_1}^2& \rho_{j_1k_2}^2  &\ldots&\rho_{j_1k_p}^2   \\
  1&  \rho_{j_2k_1}^2& \rho_{j_2k_2}^2  &\ldots&\rho_{j_2k_p}^2   \\
\vdots&\vdots &  \vdots&\ddots&\vdots \\
1&  \rho_{j_pk_1}^2& \rho_{j_pk_2}^2  &\ldots&\rho_{j_pk_p}^2
\end{array}
\right|,\\
&&
B\!
\left(
\begin{array}{ccccc}
 \star &j_1   &\ldots& j_p   \\
 0 &k_1   &\ldots&k_p   \\
\end{array}
\right) : =
\left|
\begin{array}{cccccc}
  1&r_{k_1}^2   &r_{k_2}^2&\ldots&r_{k_p}^2   \\
  1& \rho_{j_1k_1}^2& \rho_{j_1k_2}^2  &\ldots&\rho_{j_1k_p}^2   \\
  1&  \rho_{j_2k_1}^2& \rho_{j_2k_2}^2  &\ldots&\rho_{j_2k_p}^2   \\
\vdots&\vdots &  \vdots&\ddots&\vdots \\
1&  \rho_{j_pk_1}^2& \rho_{j_pk_2}^2  &\ldots&\rho_{j_pk_p}^2
\end{array}
\right|,\\
&&
B\!
\left(
\begin{array}{ccccc}
 \star&j_1   &\ldots& j_p   \\
 \star&k_1   &\ldots&k_p   \\
\end{array}
\right) : =
\left|
\begin{array}{cccccc}
  0&r_{k_1}^2   &r_{k_2}^2&\ldots&r_{k_p}^2   \\
  r_{j_1}^2& \rho_{j_1k_1}^2& \rho_{j_1k_2}^2  &\ldots&\rho_{j_1k_p}^2   \\
  r_{j_2}^2&  \rho_{j_2k_1}^2& \rho_{j_2k_2}^2  &\ldots&\rho_{j_2k_p}^2   \\
\vdots&\vdots &  \vdots&\ddots&\vdots \\
r_{j_p}^2&  \rho_{j_pk_1}^2& \rho_{j_pk_2}^2  &\ldots&\rho_{j_pk_p}^2
\end{array}
\right|,\\
&&
B\!
\left(
\begin{array}{ccccc}
 0&\star&j_1   &\ldots& j_p   \\
0 &\star&k_1   &\ldots&k_p   \\
\end{array}
\right) : =
\left|
\begin{array}{ccccccc}
0&1&1&1&\ldots&1\\
 1& 0&r_{k_1}^2   &r_{k_2}^2&\ldots&r_{k_p}^2   \\
 1& r_{j_1}^2& \rho_{j_1k_1}^2& \rho_{j_1k_2}^2  &\ldots&\rho_{j_1k_p}^2   \\
1&  r_{j_2}^2&  \rho_{j_2k_1}^2& \rho_{j_2k_2}^2  &\ldots&\rho_{j_2k_p}^2   \\
\vdots&\vdots&\vdots &  \vdots&\ddots&\vdots \\
1&r_{j_p}^2&  \rho_{j_pk_1}^2& \rho_{j_pk_2}^2  &\ldots&\rho_{j_pk_p}^2
\end{array}
\right|.
\ee
These determinants will be abbrieviated by 
$
B\!
\left(
\begin{array}{ccc}
 0 &  J    \\
0  & K    
\end{array}
\right),
B\!
\left(
\begin{array}{ccc}
 \star &  J    \\
0  & K    
\end{array}
\right),
B\!
\left(
\begin{array}{ccc}
 \star &  J   \\
\star  & K    
\end{array}
\right),
$\\
$
B\!
\left(
\begin{array}{ccccc}
 0&\star &  J    \\
0 &\star & K   
\end{array}
\right)
$
respectively.
 When $J = K,$ then 
$
B\!
\left(
\begin{array}{ccc}
 0 &  J    \\
0  & J    
\end{array}
\right), 
B\!
\left(
\begin{array}{ccc}
 \star &  J    \\
\star  & J   
\end{array}
\right), \\
$
$
B\!
\left(
\begin{array}{ccc}
 0&\star &  J    \\
0 &\star & J    
\end{array}
\right)
$ are simply written by $B(0\,J), \,B(\star J),\, B(0\star J)$
respectively.
\end{df}

For example,
$B(0 j) = -1,\, B(0jk) = 2\rho_{jk}^2,\, B(0\star j) = 2r_j^2$ and
\be
&&
B(0\star jk) = r_j^4 + r_k^4+ \rho_{jk}^4 - 2r_j^2\, r_k^2 - 2r_j^2\,\rho_{jk}^2 - 2r_k^2\,\rho_{jk}^2,\\
&&
B(0\,j\,k\,l) = \rho_{jk}^4 + \rho_{jl}^4 + \rho_{kl}^4 - 2\rho_{jk}^2\,\rho_{jl}^2
- 2\rho_{jk}^2\rho_{kl}^2 - 2\rho_{jl}^2\,\rho_{kl}^2.
\ee

From now on, we shall assume the following condition:

\medskip

\be
{ \rm({\cal H}1)}\quad(-1)^{p} \,B(0\,J) > 0,\ (-1)^{p +1}\,B(0\star J) > 0\quad(1\le p\le n+1),
\ee
for $|J| = p$ (the size of $J$ being denoted by $|J|$).

\medskip

3. Denote the inside and the outside of $S_j$  in ${\bf R}^n$ by 
\be
D_j^- : = \{f_j \le 0\},\; D_j^+ : = \{f_j \ge 0\}
\ee
respectively.
For $K \subset N,$ denote the real domains $D_K^- = \cap_{j\in K}D_j^-$ and $D_K^+ = \cap_{j\in K} D_j^+ $
respectively.

In particular, we can take as the (non-empty) domain
$D = D_{12\ldots\, n+1}^-$ defined by
\be
D_{12\ldots\,n+1}^- = \bigcap_{j=1}^{n+1}\,D_j^-, \quad
D_j^-: f_j\le0  \subset {\bf R}^n \quad(1\le j\le n+1). 
\ee

It is a non-empty spherically faced $n$-simplex 
(which will be called pseudo $n$-simplex in the sequel).
The boundary of $D$ consists of the faces $D\cap S_K$
where, for any $K \subset N$ such that $|K|=p, \,1\le p\le n,$
the intersection $S_K = \bigcap_{j\in K} S_j$ 
defines an $n-p$ dimensional sphere. In particular, 
$\bigcap_{k\in \partial_j N}\,S_k$  consists of two points.
Here the symbol $\partial_j K$ denotes the set of indices obtained by the deletion of the element $j$
from the set $K.$ $K^c$ denotes the complement of $K$ in $N.$
We can also see that $D_K^-\cap D_{K^c}^+$ are a pseudo $n$-simplex.

The orientation of ${\bf R}^n$  and $D$ is determined such that
the standard $n$ form $\varpi$ is positive:
\be
\varpi = dx_1\wedge\cdots\wedge dx_{n} > 0.
\ee

Let $P_j$ be the point in ${\bf R}^n$ such that
\be
&&\{P_j\}= \bigcap_{k\in \partial_j N} S_k\cap\partial D\quad(\partial D \;{\rm denotes\, the\,boundary\,of}\,D ).
\ee

\medskip

4. We can take the Euclidean coordinates $x_1,\ldots,x_n$
such that the polynomials $f_j$ have the following expressions:
\beq
&&f_j(x) = Q(x) + \sum_{\nu=1}^{n+1-j} \,2\alpha_{j\nu} x_\nu + \alpha_{j0}, \quad(1\le j\le n)\\
&&f_{n+1}(x) = Q(x) + \alpha_{n+1\,0}. 
\eeq
We assume for simplicity that $\alpha_{j\,n+1-j}> 0\ (1\le j\le n+1)$
and that $P_j$ satisfies the following condition:


\medskip

(${\cal H}$2)\quad The $x_{n +1- j }$-coordinate of $P_j$
is negative for every $j$.

\medskip

We have the equalities
\beq
\prod_{j=p}^{n} \alpha_{j\,n+1 - j} =  \sqrt{\frac{(-1)^{n-p} \,B(0\,p\,\ldots n\,n+1)}{2^{n-p+1}}}\quad(1\le p\le n ).
\eeq

Denote by  $\Delta[P_1,P_2,\ldots,P_{n+1}]$ and $\widetilde{\Delta}[P_1,P_2,\ldots,P_{n+1}]$
be the linear $n$-simplex and the pseudo $n$-simplex 
respectively with hyperspherical 
faces both  with vertices $P_j$ such that their sign of orientation is $(-1)^{\frac{n(n-1)}{2}}.$ 

By definition, the following properties are valid.

\bigskip

\begin{lm}

{\rm (i)} 
\be
&&(-1)^{\frac{n(n+1)}{2} + \nu -1}\,df_1\wedge\cdots \widehat{df_\nu}\cdots\wedge df_{n+1}>0\quad(1\le \nu\le n+1)
\ee
on $D$.

{\rm(ii)}
The simplex $[ P_1,P_2,\ldots,P_{n+1}]$ has the sign of orientation   
$(-1)^{\frac{n(n-1)}{2}}$ such that 
\be
[ P_1,P_2,\ldots,P_{n+1}] = (-1)^{\frac{n(n-1)}{2}} \,D.
\ee
\end{lm}

\begin{pf}
Indeed, we can show that
\beq
df_{2}\wedge\cdots\wedge df_{n+1} = 2^{n}(-1)^{\frac{(n-1)(n-2)}{2}} \prod_{j=2}^n \alpha_{j\,n+1- j}
\,x_n \,\varpi.
\eeq
Seeing that $x_n < 0,$ (i) is proved.
(ii) follows from $({\cal H}2).$
\qed
\end{pf}

\bigskip

\begin{lm}
For $J \subset N\ (|J| = p, 1\le p\le n-1),$ 
$S_J$ is a sphere of dimension $n-p$. Its radius $r_J$ equals
\be
r_J = \sqrt{- \frac{1}{2}\, \frac{B(0\star J)}{B(0\,J)}}.
\ee
\end{lm}

\bigskip

By the use of coordinates $x_j,$
$S_{n- p+2,\ldots,n+1}$ represent the spheres of dimension $n-p$
 satisfying the following equalities:
\be
&&S_{n+1} :  Q(x) = r_{n+1}^2,\\
&&S_{n - p+2\,\ldots\,n+1} : \sum_{j=p}^n\,x_j^2= r_{n - p+2\,\ldots\,n+1}^2,\\
&& S_{2\,\ldots\,n+1}  :  \{\rm two\, points\} .
\ee

\medskip

5. Denote by $f_J$ the product $\prod_{j\in J} \,f_j.$
The residues of the forms $\frac{\varpi}{f_J}$
along $S_J$ can be computed explicitly as follows.

%
\begin{prop}
For $J = \{j_1,\ldots,j_p\}\subset N,$
we have
\beq
&&{\rm Res}_{S_J}\bigl[\frac{\varpi}{f_J} \bigr] = \left[\frac{\varpi}{df_{j_1}\wedge \cdots\wedge df_{j_p}}\right]_{S_J}\nonumber \\
&&= 
\frac{(-1)^{\frac{(p-1)(p-2)}{2}}}
{\sqrt{
(-1)^{p-1}2^{p}B(0\star J)}}\varpi_J,\quad(1\le p\le n)
\eeq
where $\varpi_J$ denote the standard  sphere volume elements
on $S_J$ respectively 
such that
\beq
&&\varpi_{n+1}=  \frac{\sum_{\nu=1}^{n}\,(-1)^{\nu-1}x_\nu dx_1\wedge\cdots\widehat{dx_\nu}\cdots\wedge dx_{n}}{r_{n+1}},\nonumber\\
\\
&&\varpi_{n-p+2\,\ldots n+1} = \frac{\sum_{\nu=1}^{n-p+1} (-1)^{\nu-1}\,x_{p+\nu-1} dx_{p}\wedge\cdots\widehat{dx_{p+\nu-1}}\cdots
\wedge dx_{n}}{r_{n-p+2\,\ldots\,n+1}}, \quad(1\le p\le n-1)\nonumber\\
\\
&&\varpi_{2\,\ldots\,n+1} = -1,
\eeq
and 
\beq
r_{n-p+2\,\ldots,n+1} = \sqrt{-\frac{1}{2}\,\frac{B(0\star n-p+2\,\ldots\, n+1)}{B(0\,n-p+2\,\ldots \,n+1)}},
\eeq
$\varpi_{J}$ are obtained respectively 
from $\varpi_{n-p+2,\ldots n+1}$
by permutations of elements in  the set of indices $N.$
\end{prop}

\begin{pf}

Because of symmetry, it sufficient to prove (4) in the case where $J =\{n-p+2,\ldots, n+1\} $.

First, we prove (4) in the case of (5) and (6).
Since
\beq
df_{n-p+2}\wedge\cdots\wedge df_{n+1} = 2^{p}(-1)^{\frac{(p-1)(p-2)}{2}} \prod_{j=n-p+2}^n \alpha_{j\,n-j+1}
\sum_{j=p}^n \,x_j dx_1\wedge\cdots\wedge dx_{p-1}\wedge dx_j,\nonumber\\
\eeq
(3) (7) and (9) imply
\beq
&&df_{n-p+2}\wedge\cdots\wedge df_{n+1}\wedge \varpi_{n-p+2\ldots n+1}\nonumber\\
&&{}\quad = 2^p (-1)^{\frac{(p-1)(p-2)}{2}} \prod_{j=n-p+2}^n \alpha_{j\,n-j+1}
\frac{(\sum_{j=p}^n \,x_j^2)}{r_{n-p+2\ldots n+1}}\,\varpi\nonumber\\
&&{}\quad  = 2^p (-1)^{\frac{(p-1)(p-2)}{2}} \sqrt{\frac{(-1)^{p-1} B(0\star\,n-p+2\ldots\,n+1)}{2^{p}}}\,\varpi\quad(1\le p\le n - 1)
\nonumber.\\
\eeq

Hence, along $S_J$ it follows that 
\be
\bigl[\frac{\varpi}{df_{n-p+2}\wedge\cdots\wedge df_{n+1}}\bigr]_{S_J}
= \frac{(-1)^{\frac{(p-1)(p-2)}{2}}}
{
\sqrt{
(-1)^{p+1}\,2^p\,B(0\star\,n-p+2\ldots n+1)
}}
\, \varpi_{n-p+2\,\ldots \,n+1}.
\ee

When $p=n,$ in view of (3), (4) and $x_n < 0,$ 
we have the identity 
\be
x_n = - r_{2\ldots n+1} < 0.
\ee

Hence, at $P_1$ it follows that
\beq
\bigl[\frac{\varpi}{df_2\wedge\cdots\wedge df_{n+1}}\bigr]_{P_1}
=  -\frac{ (-1)^{\frac{(n-1)(n-2)}{2}}}{\sqrt{
(-1)^{n+1}\,2^n\,B(0\star 2\ldots n+1)}}.
\eeq
\qed
\end{pf}

\bigskip

\begin{N}
For $J \subset N,$ denote by $F_J$ the rational $n$-form
and by $W_0(J)\varpi$ a linear combination of $F_K\,(K \subset J)$
as follows:
\be
&&F_J = \frac{\varpi}{f_J},\\
&&W_0(J)\varpi = - \sum_{\nu\in J}\,B\!
\left(
\begin{array}{ccc}
  0&\star   &\partial_\nu J   \\
  0&\nu   &\partial_\nu J   \\ 
\end{array}
\right)\,F_{\partial_\nu J}
+ B(0\star J)\,F_J.
\ee
Remark that $F_J$ is also a linear combination of $W_0(K)\varpi\ (K \subset J,\,|K| \ge 1)$.
\end{N}

The following Lemma can be proved by a direct calculation (see [8] Lemma 12).

\begin{lm}

\beq
\sum_{\nu=1}^{n+1} \,(-1)^{\nu-1} \frac{df_1}{f_1}\wedge\cdots\widehat{\frac{df_\nu}{f_\nu}}\cdots
\wedge\frac{df_{n+1}}{f_{n+1}}
= \frac{2^n \,(-1)^{\frac{n(n-1)}{2}+1}}{\sqrt{(-1)^{n+1}\,2^n\, B(0 N)}}\, W_0(N)\varpi.
\eeq
\end{lm}

\bigskip

$S_{\partial_jN}$ consists of two points (denoted by $P_j,P_j^\prime$)
satisfying the equations
\be
  f_k = 0\quad(k \in \partial_jN).
\ee

The following Proposition gives the values of $f_j$ at the points $P_j, P_j^\prime.$

\begin{prop}
The values of $\frac{1}{f_j}$ at $P_j, P_j^\prime$  are negative
and positive respectively.
They are evaluated as
\beq
&&\bigl[\frac{1}{f_j}\bigr]_{P_j} = \frac{(-1)^{n+1}\sqrt{B(0\star\partial_j\!N)\,B(0 N)} + 
B
\!\left(
\begin{array}{ccccc}
  0& \star  & \partial_j N  \\
  0&j   &\partial_j N   
\end{array}
\right)
}
{B(0\star N)} < 0,\nonumber\\
\\
&&\bigl[\frac{1}{f_j}\bigr]_{P_j^\prime} = \frac{(-1)^n\sqrt{B(0\star\partial_j\!N)\,B(0 N)} + 
B
\!\left(
\begin{array}{ccccc}
  0& \star  & \partial_jN  \\
  0&j   &\partial_j
\end{array}
\right)
}
{B(0\star N)}> 0.\nonumber\\
\eeq
Due to the product formula for resultant,
\be
[\frac{1}{f_j}]_{P_j}
\,[\frac{1}{f_j}]_{P_j^\prime}
= - \frac{B(0\,\partial_jN)}{B(0\star\,N)} < 0.
\ee

\end{prop}

\begin{pf}
For simplicity, we may assume that $j = 1.$
First notice that $f_1\ne 0$ at $P_1.$
By taking the residues of both sides of (12) at $P_1$ (see [24]), we have from (11)

\beq
&&1 = {\rm Res}_{P_1} \frac{df_2}{f_2}\wedge\cdots\wedge\frac{df_{n+1}}{f_{n+1}}\nonumber\\
&&=  \frac{2^n (-1)^{\frac{n(n-1)}{2} + 1}}{\sqrt{(-1)^{n+1}\,2^n \,B(0 N)}}
\bigl\{ - B\!
\left(
\begin{array}{cccc}
  0&\star&\partial_1N     \\
 0& 1&\partial_1N    \\
\end{array}
\right)
 +
B(0\star N) \bigl[\frac{1}{f_1}\bigr]_{P_1} \bigr\}\,{\rm Res}_{P_1}
\bigl[
\frac{\varpi}{f_2\ldots f_{n+1}}\bigr]\nonumber\\
&&= \frac{ (-1)^{n + 1}}{\sqrt{B(0\star\partial_1\!N)\, \,B(0 N)}}
\bigl\{ - B\!
\left(
\begin{array}{cccc}
  0&\star&\partial_1N     \\
 0& 1&\partial_1N    \\
\end{array}
\right)
 +
B(0\star N) \bigl[\frac{1}{f_1}\bigr]_{P_1} \bigr\}.
\eeq
We can solve the equation (16) with respect to $ \bigl[\frac{1}{f_1}\bigr]_{P_1}$
and gets the formula (14).  (15) can be deduced in a similar way.
\qed
\end{pf}

\bigskip


\section{Main Theorems}

\bigskip

Main Theorems are given by \lq\lq regularization procedure of integrals"
and are a consequence from some identities  of hypergeometric integrals defined
on the $n$ dimensional complex affine space ${\bf C}^n$.

Suppose that the system of exponents $\lambda = (\lambda_1, \ldots,\lambda_{n+1)}$
are given such that all $\lambda_j > 0.$

Let $\Phi(x)$ be the multiplicative meromorphic function
\be
\Phi(x) = \prod_{j\in N} \,f_j^{\lambda_j}.
\ee

For each $J \subset N\, (1\le |J|),$  consider the integral of $|\Phi(x)|$ over the domain $D = D_J^-\cap D_{J^c}^+$:
\be
{\cal J}_\lambda (\varphi)= \int_D |\Phi(x)| \varphi \varpi,
\ee
where we take the branch of $\Phi(x)$ such that $\Phi(x) > 0$ 
for $x\in D.$
There exists a twisted $n$-cycle $\mathfrak{z}$ such that
\be
{\cal J}_\lambda(\varphi) = \int_{\mathfrak{z}} \,\Phi(x)\,\varphi\,\varpi.
\ee
  
Then the following Proposition holds true (refer to [8]) :
\begin{prop}
For each $D = D_J^-\cap D_{J^c}^+,$
the following identity holds true{\rm :}
\be
(2\lambda_\infty + n) \,{\cal J}_\lambda(1) = \sum_{p=1}^{n+1} \sum_{J \subset N, |J| = p} (-1)^p
\frac{
\prod_{j\in J} \lambda_j }{\prod_{\nu=1}^{p-1} (\lambda_\infty + n - \nu)}
\int_D |\Phi(x)| W_0(J)\,\varpi,
\ee
where  the sum ranges over the family of all unordered  sets $J$ such
that $J \subset N$ $(1 \le p\le n+1)$ and $\lambda_\infty = \sum_{j=1}^{n+1} \lambda_j$.
\end{prop}

 On the other hand, the variation of ${\cal J}_\lambda(1)$ is defined by
  \be
 d_B{\cal J}_\lambda(1) = \sum_{j=1}^{n+1}\sum_{\nu=0}^n \,d\alpha_{j\nu}\,\frac{\partial}{\partial_{\alpha_{j\nu}}}
 {\cal J}_\lambda(1).
 \ee

We want to give  an explicit variation formula for ${\cal J}_\lambda(1)$ 
with respect to the parameters $r_j^2,\, \rho_{kl}^2$. 
To do that, it is necessary to introduce the system of special 
$1$-forms $\theta_J$:

\begin{df}
\be
&&\theta_j = - \frac{1}{2}\,d\log r_j^2,\\
&&\theta_{jk}=  \frac{1}{2}\,d\log\,\rho_{jk}^2,\\
&&\theta_J = (-1)^p \, \sum_{j,k\in J, j<k} \frac{1}{2}\, d\log B(0 j k)\cdot\\
&&\sum_{\mu_1,\ldots,\mu_{p-2}}
\prod_{\nu=1}^{p-2}
\frac{
B\!
\left(
\begin{array}{cccccccc}
 0 & \star  &\mu_{\nu-1}&\ldots& \mu_1&j&k   \\
 0 &  \mu_\nu & \mu_{\nu-1}&\ldots&\mu_1&j&k  \\
\end{array}
\right)
}
{B(0\mu_\nu\mu_{\nu-1}\ldots\mu_1\,j\,k)},\quad(2\le p\le n+1,\, |J| = p)
\ee
where $\mu_1,\ldots,\mu_{p-2}$ ranges over the family of all ordered sequences
consisting of $p-2$ different elements of $\partial_{j}\partial_k J.$
\end{df}

Then we have the following (refer to  [8]).
\begin{prop}
For each $D = D_J^-\cap D_{J^c}^+,$
we have
\be
d_B{\cal J}_\lambda(1)  = \sum_{p=1}^{n+1} \sum_{J\subset N, |J|=p}\,
\frac{\prod_{j\in J}\lambda_j}{\prod_{\nu=1}^{p-1}(\lambda_\infty + n - \nu)}\,\theta_J\,
\int_D \,|\Phi(x)| W_0(J)\varpi.
\ee
$J$ ranges over the family of unordered sets such that $J \subset N$
and $1 \le p \le n+1.$
\end{prop}
 
\medskip

We now take $D = D_N^-.$
The following Theorem is an immediate consequence of  Propositions 8 and 9
tending  $\lambda_j \to 0$ for all positive $\lambda_j.$

\bigskip

\bigskip

[{\bf  Theorem I}]

\medskip

{\it
Assume the conditions $({\cal H }1)$ and $({\cal H}2)$. 
Let $v(D)$ be the volume of 
the pseudo $n$-simplex 
$D = D_{N}^-:$
\be
v(D) = \int_D \varpi .
\ee

Then

{\rm (i)}
we have the identity\text{:}
\beq
&&n\, v(D) = - \sum_{p=1}^n 
\frac{(n-p)!}{(n-1)!} 
\sum_{J\subset N, |J| = p}
(-1)^p\,\sqrt{\frac{(-1)^{p+1}\,B(0\star J)}
{2^p} }\,v_J \nonumber\\
&&{}\qquad \qquad + \,(-1)^{n}\,  \frac{1}{(n-1)!}\, \sqrt{\frac{(-1)^{n+1} B(0 N)}{2^{n}}}.
\eeq

{\rm(ii)}
the following variational formula holds true:
\beq
&&d_B v(D) = - \sum_{p=1}^n 
\frac{(n-p)!}{(n-1)!} 
\sum_{J\subset N, |J| = p}
\sqrt{\frac{(-1)^{p+1}\,B(0\star J)}
{2^p} }\,\theta_J\,v_J \nonumber\\
&&{}\qquad \qquad -  \frac{1}{(n-1)!}\, \sqrt{\frac{(-1)^{n+1} B(0 N)}{2^{n}}}\; \theta_N.
\eeq

Here $v_j,\,v_{jk},\, v_J$ denote the lower dimensional  volumes corresponding to
the boundaries $S_j\cap \partial D, S_j\cap S_k\cap \partial D,\, S_J\cap \partial D\,(|J| \ge 3)$
respectively  such that 
\beq
 v_J =  \int_{S_J\cap \partial D} |\varpi_J|\quad(1\le |J| \le n).
\eeq
\it}

\begin{re}
The formula {\rm (18)} is just an analog of the classical variational formula due to L.Schl\"afli
concerning the volume of a $n$ dimensional geodesic simplex in the unit hypersphere.

\end{re}

\begin{pf}

Since both identities can be proved in the same way, we only give a proof for the latter
identity (18).

{\rm Proposition 9} shows
\beq
&&d_B {\cal J}_\lambda(1) =  \sum_{j=1}^{n+1}\,\lambda_j  \int_D |\Phi(x)| W_0(j)\varpi\,\theta_j
+ \sum_{j<k}\frac{\lambda_j\lambda_k}{\lambda_\infty + n - 1} \int_D|\Phi(x)| W_0(jk)\varpi\,\theta_{jk}\nonumber\\
&&{}\qquad \qquad + \sum_{p=3}^n 
\frac{\prod_{j\in J}\,\lambda_j}
{\prod_{\nu=1}^{p-1}\,(\lambda_\infty + n - \nu)}
 \int_D\,|\Phi(x)| 
 W_0(J)\varpi\,\theta_{J}\nonumber\\
&&{}\qquad \qquad +  \frac{\prod_{j\in N}\,\lambda_j}
{\prod_{\nu=1}^n(\lambda_\infty+n - \nu)}
\int_D\, |\Phi(x)|\,W_0(N)\varpi\,\theta_N.
\eeq

Let us take the limit for $\lambda_j \to 0 \ (1\le j\le n+1)$
 on both sides of (19) In the {\rm LHS} such that 
\be
\lambda_j = \varepsilon \,(\varepsilon \downarrow 0).
\ee

In the {\rm LHS},
\be
\lim_{\lambda\to 0\,(1\le j\le n+1)}d_B {\cal J}_\lambda(1) = d_Bv(D).
\ee

In the {\rm RHS}, first remark that
\be
\lim_{\varepsilon \downarrow  0} \prod_{j\in J}\,\lambda_j \int_D |\Phi(x)| \varphi(x) F_K = 0,
\ee
provided $K \varsubsetneq J.$

In  the {\rm RHS}, seeing that $f_j < 0$ in $D$,
due to {\rm Proposition 4}, the following equalities hold by the method of generalized functions
(see [11] Chap III, 2):
\beq
&&\lim_{\lambda\to 0}\lambda_j \int_D |\Phi(x)|\,\varphi(x)\,W_0(j)\,\varpi = 
\lim_{\lambda\to 0}\lambda_j \int_D |\Phi(x)|\,\varphi(x)\,
B(0\star j)\,\frac{\varpi}{f_j}\nonumber\\
&&{}\qquad \qquad= - B(0\star j)\,\lim_{\lambda\to 0}\lambda_j \int_D |\Phi(x)|\,\varphi(x)
\frac{\varpi}{|f_j|} \nonumber\\
&&{}\qquad \qquad =
- \sqrt{\frac{B(0\star j)}{2}}\,\int_{S_j\cap \partial D}[\varphi]_{S_j}\,|\varpi_j|,\\
&&\lim_{\lambda\to 0}\lambda_j \lambda_k
\int_D |\Phi(x)| \,
\varphi(x)\,W_0(jk)\,\varpi
= B(0\star jk)\,\lim_{\lambda\to 0}\lambda_j \lambda_k
\int_D |\Phi(x)| \,
\varphi(x)\,\frac{\varpi}{f_jf_k} \nonumber\\
&&{}\qquad \qquad = B(0\star jk)\,\lim_{\lambda\to 0}\lambda_j \lambda_k
\int_D |\Phi(x)| \,
\varphi(x)\,\frac{\varpi}{|f_jf_k|} \nonumber\\
&&{}\qquad \qquad =B(0\star jk)\,
\int_D |\Phi(x)| \,
\varphi(x)\,|\frac{\varpi}{df_j\wedge df_k}| \nonumber\\
&&{}\qquad \qquad =
- \sqrt{\frac{- B(0\star jk)}{4}}\int_{S_j\cap S_k\cap \partial D}[\varphi]_{S_{jk}}\,
|\varpi_{jk}|.
\eeq
In general, we have, for $|J| \le n,$
\beq
&&\lim_{\lambda\to 0}\prod_{j\in J}\,\lambda_j
\int_D 
|\Phi(x)|\,\varphi(x)\,W_0(J)\varpi
 = B(0\star J) \lim_{\lambda\to 0}\prod_{j\in J}\,\lambda_j
\,\int_D 
|\Phi(x)|\,\varphi(x)\,
F_{J}\nonumber\\
&&{}\qquad \qquad = (-1)^p\,B(0\star J)
\lim_{\lambda\to 0}\prod_{j\in J}\,\lambda_j
\int_D 
|\Phi(x)|\,\varphi(x)\,
|F_{J}|\nonumber\\
&&{}\qquad \qquad =
(-1)^p\,B(0\star J)\,\int_{S_J\cap\partial D}\,[\varphi]_{S_J}
\bigl|\frac{
\varpi}{
df_{j_1}\wedge \cdots\wedge df_{j_p}}\bigr|_{S_J}  \nonumber\\
&&{}\qquad \qquad = - \sqrt{\frac{(-1)^{p+1}\,B(0\star J)}{2^p}}
\int_{S_J\cap\partial D}\,[\varphi]_{S_J} [\varphi]_{S_J} |\varpi_J|,
\eeq 
where $J = \{j_1,\ldots,j_p\},\,p \le n.$

As for the last  term of the {\rm RHS} of {\rm(20)}, 
when $\lambda\to 0,$ the limit value is divided into  the ones at $P_j\ (1\le j\le n+1).$

We assume that $D$ is divided into $n+1$ domains $D_j^*\,(j\in  N)$ such that
$P_j$ lies only in the inside  of  $\partial D_j^*$:
\be
D = \bigcup_{j\in N}\,D_j^*.
\ee
Then
\be
\int_D |\Phi(x)| \varpi = \sum_{j\in N} \,\int_{D_j^*} |\Phi(x)|\,\varpi.
\ee

 We restrict ourselves to  the integral over $D_1^*.$ 
 Since $f_k < 0\ (1\le k\le n+1)$ in the inside of $D_1^*$
and $\bigl[f_k\bigr]_{P_1} = 0 \ (2\le k\le n+1)$  
and $\bigl[f_1\bigr]_{P_1} < 0 ,$
\be
&&\lim_{\lambda \to 0}
\frac{\prod_{j=1}^{n+1}\,\lambda_j}{\prod_{\nu=0}^{n-1} (\lambda_\infty + \nu)}
\int_{D_1^*} |\Phi(x)|\, \,W_0(N)\,\bigl|\varpi\bigr| \\
&&{}\qquad  =
\lim_{\varepsilon \downarrow 0}
 \frac{(-1)^n}{(n+1)\,\prod_{\nu=1}^{n-1}((n+1)\varepsilon + \nu)}\cdot \\
&& {}\quad \quad \left(
- B\!
\left(
\begin{array}{ccccc}
 0 &\star   & \partial_1 N\\
 0&1&\partial_1 N
\end{array}
\right) +  B(0\star N) [\frac{1}{f_1}]_{P_1}\right) 
|\bigl[\frac{\varpi}{df_2\wedge \cdots\wedge df_{n+1}}\bigr]_{P_1}|\\
&&{}\qquad  = -\,\frac{1}{(n+1)\,(n-1)!} \sqrt{\frac{(-1)^{n+1}\,B(0 N)}{2^n}},
\ee
in view of (12) and Proposition 6.

Similarly, the integral over each $D_j^*$ has the same value as above:
\beq
&&\lim_{\varepsilon \downarrow0}
\frac{\prod_{j=1}^{n+1}\,\lambda_j}{\prod_{\nu=0}^{n-1} (\lambda + \nu)}
\int_{D_j^*} |\Phi(x)|\, \,W_0(N)\,\varpi \nonumber\\
&&{}\qquad = - \, \frac{1}{(n+1)\,(n-1)!} \sqrt{\frac{(-1)^{n+1}\,B(0 N)}{2^n}}.
\eeq

(20) and (24) imply Main Theorem {\rm I (ii)}.
{\rm (i)} can be proved in the same way from Proposition 7. 
\qed
\end{pf}

\medskip

For $D = D_J^-\cap D_{J^c}^+$, the similar formulae to (17) and (18)
for $v(D)$ can be described and be proved in the same way:
\be
&&nv(D) = - \sum_{p=1}^n \frac{(n-p)!}{(n-1)!}\, \sum_{J \subset N, |J| = p}
(-1)^p \sqrt{\frac{(-1)^{p+1}\,B(0\star J)}{2^p}} \,v_J \\
&&{}\qquad \qquad + (-1)^n \frac{1}{(n-1)!} \sqrt{\frac{(-1)^{n+1} B(0\,N)}{2^n}}.
\ee
\be
&&d_B\,v(D) = - \sum_{p=1}^n \frac{(n-p)!}{(n-1)!}\, \sum_{J \subset N, |J| = p}
(-1)^p \sqrt{\frac{(-1)^{p+1}\,B(0\star J)}{2^p}} \,\theta_J\,v_J \\
&&{}\qquad \qquad -  \frac{1}{(n-1)!} \sqrt{\frac{(-1)^{n+1} B(0\,N)}{2^n}}\,\theta_N.
\ee

\bigskip

In place of $({\cal H}1),$ we can also consider the following condition:

\be
&&({\cal H}1') \quad (-1)^p\, B(0\,J) > 0,\; (-1)^{p+1}\,B(0\star J) > 0\quad(1\le p\le n\, \mbox{ for}\, |J| = p),\\
&&\hspace{1.5cm}  (-1)^{n+1}\,B(0\,N) > 0,\, (-1)^{n} \,B(0\star N) < 0,
\ee
instead of $(-1)^n\, B(0\star N) > 0$ (For $n=2,$ refer to Figure 2).

\bigskip

There are $2^{n+1} - 1$ non-empty bounded chambers 
\be
D_J^- :  f_j \le 0 \quad(j\in J),\ f_j \ge 0 \quad(j\in J^c),
 \ee
where $J$ ranges over the family of subsets in $N$
such that $|J| = p, 0 \le p \le n  $.
When $p = 0,$ i.e.,  when $J$ is empty, we may denote $D_J^-$ by $D_N^+$:
\be
D_N^+ :  f_j \ge 0\quad(1\le j\le n+1).
\ee
We assume that the inside of $D_N^+$ is a non-empty domain.
$D_N^+$ is the support of the $n$-simplex $\tilde{\Delta}[P_1\,\ldots\,P_{n+1}]$
such that $P_j \in S_{\partial_j N}.$
$S_{\partial_j N }$ consists of  two points $\{P_j,\,P_j^\prime\}$
(see Figure 2 for 2 dimensional case ). Then
\be
[\frac{1}{f_j}]_{P_j} > 0, \,[\frac{1}{f_j}]_{P_j^\prime} > 0,
\ee
and
\be
\bigl[\frac{1}{f_j}\bigr]_{P_j}\,\bigl[\frac{1}{f_j}\bigr]_{P_j^\prime} = - \frac{B(0\,\partial_jN)}{B(0\star N)} > 0.
\ee

For the volume of $D_N^+,$
the following theorem holds true.

\bigskip

\bigskip

[{\bf Theorem II}]

\medskip

{\it
Let $v = v(D_N^+)$ be the volume of $D_N^+$
and $v_J$ be similarly defined as in {\rm(19)} by
\be
&&v = v(D_N^+) = \int_{D_N^+} \varpi,\\
&&v_J = \int_{S_J\cap D_N^+}\, |\varpi_{J}|\quad(1\le p\le n,\, p = |J|).
\ee

Then the following identities hold true{\rm :}

{\rm(i)}
\be
&&n\,v(D_N^+) =  - \sum_{p=1}^n \frac{(n-p)!}{(n-1)!}\,\sum_{|J| = p}  \sqrt{\frac{(-1)^{p+1}\,B(0\star J)}{2^p}}\,
v_J \\
&&{}\qquad \qquad  + \frac{1}{(n-1)!} \,\sqrt{\frac{
 (-1)^{n+1} \,B(0\,N)}{2^n}
 }.
 \ee
 
 {\rm(ii)} 
 \be
&& d_Bv(D_N^+) = - \sum_{p=1}^n \frac{(n-p)!}{(n-1)!}\,\sum_{|J| = p} (-1)^p\,\theta_J\, \sqrt{\frac{(-1)^{p+1}\,B(0\star J)}{2^p}}\,
v_J \\
&&{}\qquad \qquad  +\, \frac{(-1)^{n+1}}{(n-1)!} \,\sqrt{\frac{
 (-1)^{n+1} \,B(0\,N)}{2^n}}\,\theta_N
 .
 \ee
}

\begin{pf}
Theorem II can also be  proved in the same way as Theorem I. 
This theorem also follows from the identities in Theorem 1 by an analytic continuation 
(i.e., Picard-Lefschetz transformation
of twisted cycles ) of $v(D_N^-)$
with respect to the parameters $r_j^2,\,\rho_{jk}^2$ such that
\be
&& B(0\star N) \longrightarrow \,- B(0\star N),\\
&&v(D_N^-) \longrightarrow (-1)^n\, v(D_N^+),\\
&&v(D_J^-) \longrightarrow (-1)^{p-1}\, v(D_J^-).
\ee
\qed
\end{pf}

\begin{re}
An elementary proof of Theorem II (i) (and therefore of Theorem I (i)) will be given in the Appendix.
\end{re}

In the following, we give a simple application of Main Theorems.

\bigskip

\section{Examples}


\bigskip

[{\bf example1}]

\medskip

In the $n$ dimensional Euclidean space, consider two hyperspheres $S_1,\,S_2$ with the centers
$O_1,\,O_2$  and with radii $r_1,\,r_2$ such that the distance between $O_1$ and $O_2$
is equal to $\rho_{12}.$

Assume $S_1\cap S_2$ is a non-empty $n-2$ dimensional sphere. 
$S_1\cap S_2$ is contained in the hyperplane $L$ which intersects the segment $\overline{O_1O_2}$
at a point $M.$ 

The radius $h$ of $S_1\cap S_2$, the distance $O_1M$ and $O_2M$ are expressed as
\be
&&h = \frac{\sqrt{- B(0\star 12)}}{2\rho_{12}} = r_1\sin \frac{1}{2}\psi_{12}  = r_2\sin\frac{1}{2}\psi_{21},\\
&& O_1M = r_1\cos\frac{1}{2}\psi_{12}= \frac{B\!
\left(
\begin{array}{ccc}
 0 & 2  &1   \\
 0 &\star   &1   \\
\end{array}
\right)
}{2\rho_{12}}, \\
&&O_2M =  r_2\cos\frac{1}{2}\psi_{21}
=\frac{B\!
\left(
\begin{array}{ccc}
 0 & 1  &2   \\
 0 &\star   &2   \\
\end{array}
\right)
}{2\rho_{12}}
\ee
such that
\be
 \rho_{12} = r_1\cos\frac{1}{2}\psi_{12}+ r_2\cos\frac{1}{2}\psi_{21},
 \ee
where $\psi_{12},\psi_{21}$ satisfy $0 < \psi_{12} < \pi,\, 0<\psi_{21}<\pi.$

Denote by $D_{12}^-$ the common domain (lens domain) surrounded by $S_1, \,S_2$.
$S_1\cap S_2$ is an $(n-2)$ dimensional sphere. 

The volume $v(D_{12}^-)$ of $D_{12}^-$ can be evaluated by an elementary calculus
as follows:
\beq
v(D_{12}^-) = v_1 + v_2,
\eeq
where
\be
&& v_1 =\frac{1}{n-1}\, C_{n-2} r_1^n \, \int_{\cos\frac{1}{2}\psi_{12}}^1 (1- \tau^2)^{\frac{n-1}{2}}\,d\tau,\\
&&v_2 =  \frac{1}{n-1}\,C_{n-2} r_2^n \, \int_{\cos\frac{1}{2}\psi_{21}}^1 (1- \tau^2)^{\frac{n-1}{2}}\,d\tau,
\ee
$v_1,v_2$ denote the volumes of the domains surrounded by $S_1, L$
and $S_2, L$ respectively, and $C_{n-2}$ denotes
the volume of the $n-2$ dimensional unit hypersphere:
\be
C_{n-2} =\frac{
 2\pi^{
\frac{n-1}{2}}
}{\Gamma(\frac{n-1}{2})}.
\ee

The lower dimensional volumes of $S_1\cap S_2,\, S_1\cap \partial D_{12}^-,\, S_2\cap \partial D_{12}^-$
equal respectively
\beq
&&v(S_1\cap S_2) = C_{n-2}\,h^{n-2},\nonumber\\
&&v(S_1\cap\partial D_{12}^-) =\frac{\partial v}{\partial r_1}
= C_{n-2}\,r_1^{n-1} \int_0^{\frac{\psi_{12}}{2}} \,\sin^{n-2}t\,dt ,\\
&&v(S_2\cap\partial D_{12}^-) =\frac{\partial v}{\partial r_2}
= C_{n-2}\,r_2^{n-1} \int_0^{\frac{\psi_{21}}{2}} \,\sin^{n-2}t\,dt.
\eeq
The integral in the {\rm RHS} can be expressed in the 
following expansion:
\beq
&&\int_0^{\frac{\psi_{jk}}{2}} \,\sin^{n-2}t\,dt
 = - \sum_{0\le 2\nu\le n-3} \cos\frac{\psi_{jk}}{2}\,\frac{(n-3)\cdots(n-2\nu+1)}{(n-2)\cdots(n-2\nu)}
\bigl(\sin\frac{\psi_{jk}}{2}\bigr)^{n-3-2\nu}\nonumber\\
&& {}\qquad \qquad \qquad + 
\left\{
\begin{array}{ccc}
  C_{n-2}^\prime (1- \cos\frac{\psi_{jk}}{2})\\
  C_{n-2}^\prime \frac{\psi_{jk}}{2}
\end{array}
\right.\quad(\{j,k\} = \{1, 2\}),
\eeq
where $2C_{n-2}^\prime$ or $2\pi\,C_{n-2}^\prime$ equals
$\frac{\sqrt{\pi}\,\Gamma(\frac{n-1}{2})}{\Gamma(\frac{n}{2})}$
according as $n$ is odd or even.

\begin{re}
In the case where $n$ is odd, $\dim S_j\cap \partial D_{12}^-\,(j = 1,2)$ is even, 
{\rm (25) - (27)} are related with the generalized Gauss-Bonnet formula.
Indeed, 
the second formula due to Allendoerfer-Weil (see [2]) can be applied
to $v(S_1\cap \partial D_{12}^-)$ or $v(S_2\cap \partial D_{12}^-)$.  
The above formulae coincide with it.

The derivation of $v(D_{12}^-)$ with respect to $r_1,r_2,\rho_{12}$ in  (25) leads to
 the following formula
 \beq
&& dv(D_{12}^-) = v(S_1\cap\partial D_{12}^-) \,dr_1  + v(S_2\cap\partial D_{12}^-) \,dr_2 \nonumber\\
&&{}\qquad \qquad  - \frac{1}{n-1}\, v(S_1\cap S_2) 
 \sqrt{\frac{- B(0\star12)}{4}}\,\theta_{12}\nonumber\\
 \eeq
 with $\theta_{12} = \frac{1}{2}\,\log \rho_{12}^2.$

(29) is nothing else than a special case  in the $n$ dimensional space 
derived from (18) for $\varepsilon\to 0,$
after putting to be $\lambda_1 = \lambda_2 = \varepsilon$
and 
 $\lambda_j = 0\ (3\le j\le  n+1)$.
\end{re}

\medskip

In particular, in the case where $n = 2, 3,$ the volumes $v(D_{12}^-)$ are simply written as
\beq
&&\bullet \;v(D_{12}^-) = \frac{1}{2}\, r_1^2\,(\psi_{12} - \sin \psi_{12}) + \frac{1}{2}\,
\,r_2^2 \,(\psi_{21} - \sin \psi_{21})\quad( n = 2),\\
&&\bullet\;v(D_{12}^-) = \pi r_1^3 \,(\frac{2}{3} - \cos\frac{1}{2}\psi_{12} + \frac{1}{3}\,\cos^3\frac{1}{2}\,\psi_{12})\nonumber\\
&&{}\qquad \qquad + \pi r_2^3 \,(\frac{2}{3} - \cos\frac{1}{2}\psi_{21} + \frac{1}{3}\,\cos^3\frac{1}{2}\psi_{21})\quad( n = 3),\nonumber\\
\eeq
where
$\psi_{12}, \psi_{21}$ denote the angles at $O_1, O_2$ respectively 
subtended by the diameter
of $S_1\cap S_2.$
Remark that 
\be
&&e^{i\frac{\psi_{12}}{2}}  = \frac{B\!
\left(
\begin{array}{ccc}
 0 & \star  &1   \\
 0 &2   &1   \\
\end{array}
\right)
+ i\sqrt{-B(0\star 12)}}{2\rho_{12}\,r_1} ,\\
&&e^{i\frac{\psi_{21}}{2}}  = \frac{B\!
\left(
\begin{array}{ccc}
 0 & \star  &2   \\
 0 &1   &2   \\
\end{array}
\right)
+ i\sqrt{-B(0\star 12)}}{2\rho_{12}\,r_2} .
\ee

In the case where $n =2, 3,$ the formula (28) becomes 
\beq
&&\bullet\;dv(D_{12}^-) = r_1\psi_{12}\,dr_1 + r_2\psi_{21}\,dr_2 - \sqrt{- B(0\star12)}\, \frac{d\rho_{12}}{\rho_{12}},\\
&&\bullet\;dv(D_{12}^-) = \frac{\pi r_1}{\rho_{12}}\{r_2^2 - (r_1- \rho_{12})^2\}\,dr_1
+ \frac{\pi r_2}{\rho_{12}}\{r_1^2 - (r_2- \rho_{12})^2\}\,dr_2\nonumber\\
&&{}\qquad \qquad  - \frac{\pi}{4\rho_{12}^2} \,B(0\star 12)\,d\rho_{12},\nonumber\\
\eeq
in view of the identity
\beq
\frac{1}{2}\,d\psi_{jk} = \frac{1}{\sqrt{- B(0\star jk)}}
\{- B\!
\left(
\begin{array}{ccc}
 0 &\star   & j  \\
  0&\star   &k   
\end{array}
\right)\,\frac{dr_j}{r_j} + 2r_k\,dr_k - 
B\!
\left(
\begin{array}{ccc}
 0 &  j & k  \\
 0 &\star   & k  
\end{array}
\right)
\,\frac{d\rho_{jk}}{\rho_{jk}}\}\nonumber\\
\eeq
for $j,k = 1,2$ or $2,1$ respectively.


\bigskip

\bigskip

[{\bf Example 2}]

\medskip

Assume that $n=2$.  Then $D_{123}^-$ is the pseudo-triangle $\tilde{\Delta}[P_1P_2P_3]$
with vertices
$P_1= (\xi_1,\xi_2)),\,P_2(\eta_1,\eta_2),\,P_3(\zeta_1,\zeta_2)$ (see Figure 1),
where
\be
&&\xi_1 = -
 \frac{B\!
\left(
\begin{array}{ccc}
 0 &2   &  3 \\
  0&\star   &3   \\ 
\end{array}
\right)
}{\sqrt{2B(023)}},\,
\xi_2 = - \sqrt{\frac{- B(0\star 23)}
{2B(023)}},\\
&&\eta_1 = \frac{1}{B(013)\sqrt{2B(023)}}
\bigl\{- B\!
\left(
\begin{array}{ccc}
 0 &1   &3   \\
  0&\star   &3   \\
\end{array}
\right)
B\!
\left(
\begin{array}{ccc}
 0 &1   &3   \\
  0&2   &3   \\
\end{array}
\right)
- \sqrt{B(0\star13)\,B(0123)}\bigr\},\\
&&\eta_2 = \frac{1}{B(013)\sqrt{2B(023)}}
\{- B\!
\left(
\begin{array}{ccc}
 0 &1   &3   \\
  0&\star   &3   \\
\end{array}
\right)\sqrt{- B(0123)}\}
- B\!
\left(
\begin{array}{ccc}
 0 &1   &3   \\
  0&2   &3   \\
\end{array}
\right)
\sqrt{- B(0\star13)}\},\\
&&\zeta_1 = 
 \frac{1}{B(012)\sqrt{2B(023)}}
\{
B\!
\left(
\begin{array}{ccc}
 0 &1   &2   \\
  0&\star   &2   \\
\end{array}
\right)
B\!
\left(
\begin{array}{ccc}
 0 &1   &2   \\
  0&3   &2   \\
\end{array}
\right)
- B(012)\,B(023)\\
&&{}\qquad  + \sqrt{B(0123)\,B(0\star12)}\},\\
&&\zeta_2 =  \frac{1}{B(012)\sqrt{2B(023)}}
\{- B\!
\left(
\begin{array}{ccc}
 0 &1   &2   \\
  0&\star   &2   \\
\end{array}
\right)\sqrt{- B(0123)}
+ B\!
\left(
\begin{array}{ccc}
 0 &1   &2   \\
  0&3   &2   \\
\end{array}
\right)
\sqrt{- B(0\star12)}\}.
\ee
Note that $\xi_2 < 0,\, \eta_1 < 0.$

The area of $\Delta(O_1O_3O_2)$ is expressed by
\beq
|\Delta(O_1O_3O_2)| = \frac{1}{2}\, |\delta|,
\eeq 
where $\delta$ denotes
\be
\delta = \left|
\begin{array}{ccc}
1  &\xi_1   &\xi_2   \\
 1 & \eta_1  &\eta_2   \\
 1 & \zeta_1  & \zeta_2  
\end{array}
\right| = 
- \frac{1}{2} \,\sqrt{- B(0123)}< 0.
\ee

Denote by $\varphi_j$ the angle of the triangle $\Delta(O_1O_3O_2)$
at the vertex $O_j$.
Then 
\beq
e^{i\varphi_j} = \frac{B\!
\left(
\begin{array}{ccc}
0  & k  &j   \\
 0 &l   & j  \\
\end{array}
\right) 
+ i\sqrt{- B(0123)}
}{2\rho_{jk}\rho_{jl}}\quad(j,k,l \,\mbox{different\, indices}).
\eeq

Denote by $P_1^\prime,\, P_2^\prime, \,P_3^\prime$
the intersection points of $S_2\cap S_3, \,S_3\cap S_1,\,
S_1\cap S_2$ which are different from $P_1,\,P_2,\,P_3$ respectively.
Also denote by $\psi_{jk}$ the angle at $O_j$ subtended by the arc $\widehat{P_k P_k^\prime}\cap S_j$.

Then Theorem I (i) shows

\begin{lm}
\beq
&&v(D) = \tilde{\Delta}(P_1P_3P_2)\nonumber\\
&& = |\Delta(O_1O_3O_2)| - \sum_{j=1}^ 3\, |\Delta(O_1O_3O_2)\cap D_j |+ \sum_{1\le j<k\le 3} 
|\Delta(O_1O_3O_2)\cap D_j^-\cap D_k^-|, \nonumber\\
\eeq
where owing to {\rm(30)}
\be
&&|\Delta(O_1O_3O_2)\cap D_j^-| = \,\frac{1}{2}\,r_j^2\varphi_j,\\
&& |\Delta(O_1O_3O_2)\cap D_j^-\cap D_k^- |= \frac{1}{2}\,| D_j^-\cap D_k^-|\\
&&{}\qquad  =\frac{1}{4} \,r_j^2\,(\psi_{jk} - \sin\psi_{jk}) + \frac{1}{4} \,r_k^2\,(\psi_{kj} - \sin\psi_{kj}).
\ee
\end{lm}

Denote by $\psi_1,\,\psi_2,\,\psi_3$ the angles at $O_j$ subtended 
by the sides $\widehat{P_2P_3}, \widehat{P_3P_1},\widehat{P_1P_2}$
of  the pseudo triangle $\tilde{\Delta}(P_1,P_3,P_2)$
respectively such that the arc length $s_j$ of $\widehat{P_k P_l}$
is equal to
\be
s_j = r_j\,\psi_j.
\ee

 $\psi_j$ are also related with $\psi_{jk}, \,\varphi_j$
as follows:
\be
\psi_j = \frac{1}{2}\,\psi_{jk} + \frac{1}{2}\,\psi_{jl} - \varphi_j.
\ee

The formula (37) can be identified with the special case $n = 2$ of  (17).

On the other,  since $\varphi_1 + \varphi_2 + \varphi_3 = 2\pi$,
\be
&&\psi_1 + \psi_2 + \psi_3 = \frac{1}{2} \,\sum_{j\ne k} \psi_{jk}\\
&& {}\qquad = 2\pi - \angle P_1P_3P_2 - \angle P_2P_1P_3 - \angle P_3P_2P_3.
\ee
This identity is a special case of the second Allendoerfer-Weil formula 
in the Euclidean plane (see [2] Theorem II).
Furthermore,
\be
2v(D) = r_1\,v_1 + r_2\,v_2 + r_3\,v_3 -  \frac{1}{2}\,\sum_{1\le j<k\le 3}\, \sqrt{- B(0\star\,j\,k)}\,v_{jk}
+ \frac{1}{2} \sqrt{- B(0\,1\,2\,3)},
\ee
with $v_j = r_j \psi_j$. This identity coincides with (17) in the two dimensional case.

Taking into consideration the identities (30),  (31), (33) and the following equalities
($ j,k,l\, {\rm are \,different\, indices\,of\,  } 1,2,3$)
\be
&&dB(0123) =  \sum_{j<k}\,
\frac{1}{2}\, d\rho_{jk}\, B\!
\left(
\begin{array}{ccc}
 0 &j   & l  \\
 0 &k  &l   
\end{array}
\right),\\
&&d\varphi_j = \frac{1}{\sqrt{- B(0123)}}\,
\{- B\!
\left(
\begin{array}{ccc}
  0&j   &k  \\
  0&l   &k   \\ 
\end{array}
\right)\,\frac{d\rho_{jk}}{\rho_{jk}}
- B\!
\left(
\begin{array}{ccc}
  0&j   &l  \\
  0&k   &l   \\ 
\end{array}
\right)\,\frac{d\rho_{jl}}{\rho_{jl}} + 2\rho_{kl}d\rho_{kl}\},
\ee
we get the formula 
\be
&&dv(D) =  \sum_{j=1}^3 r_j\,\psi_j \, dr_j 
- \frac{1}{2} \sum_{j<k}\sqrt{- B(0\star jk)}\,\frac{d\rho_{jk}}{\rho_{jk}}\\
&&{}\qquad - \,\frac{1}{2\sqrt{-B(0123)}}\,
\{\sum_{j<k}^3 B\!
\left(
\begin{array}{ccccc}
 0 &\star   &j&k  \\
 0 & l  &j&k   \\
\end{array}
\right)
\,\frac{d\rho_{jk}}{\rho_{jk}}\},
\ee
which is nothing else than (18) for $n= 2$
in view of the identity
\be
&&\theta_{123} = - \frac{1}{B(0123)}\,
\{B\!
\left(
\begin{array}{cccc}
 0 &\star   &1&2   \\
 0 &3   &1&2   
\end{array}
\right)
\frac{d\rho_{12}}{\rho_{12}} +
B\!
\left(
\begin{array}{cccc}
 0 &\star   &1&3  \\
 0 &2   &1&3   
\end{array}
\right)
\frac{d\rho_{13}}{\rho_{13}}\\
&&{}\qquad +
B\!
\left(
\begin{array}{cccc}
 0 &\star   &2&3   \\
 0 &1   &2&3   
\end{array}
\right)
\}\frac{d\rho_{23}}{\rho_{23}}.
\ee
\begin{figure}[htbp]
 \begin{center}
  \includegraphics[width=390pt]{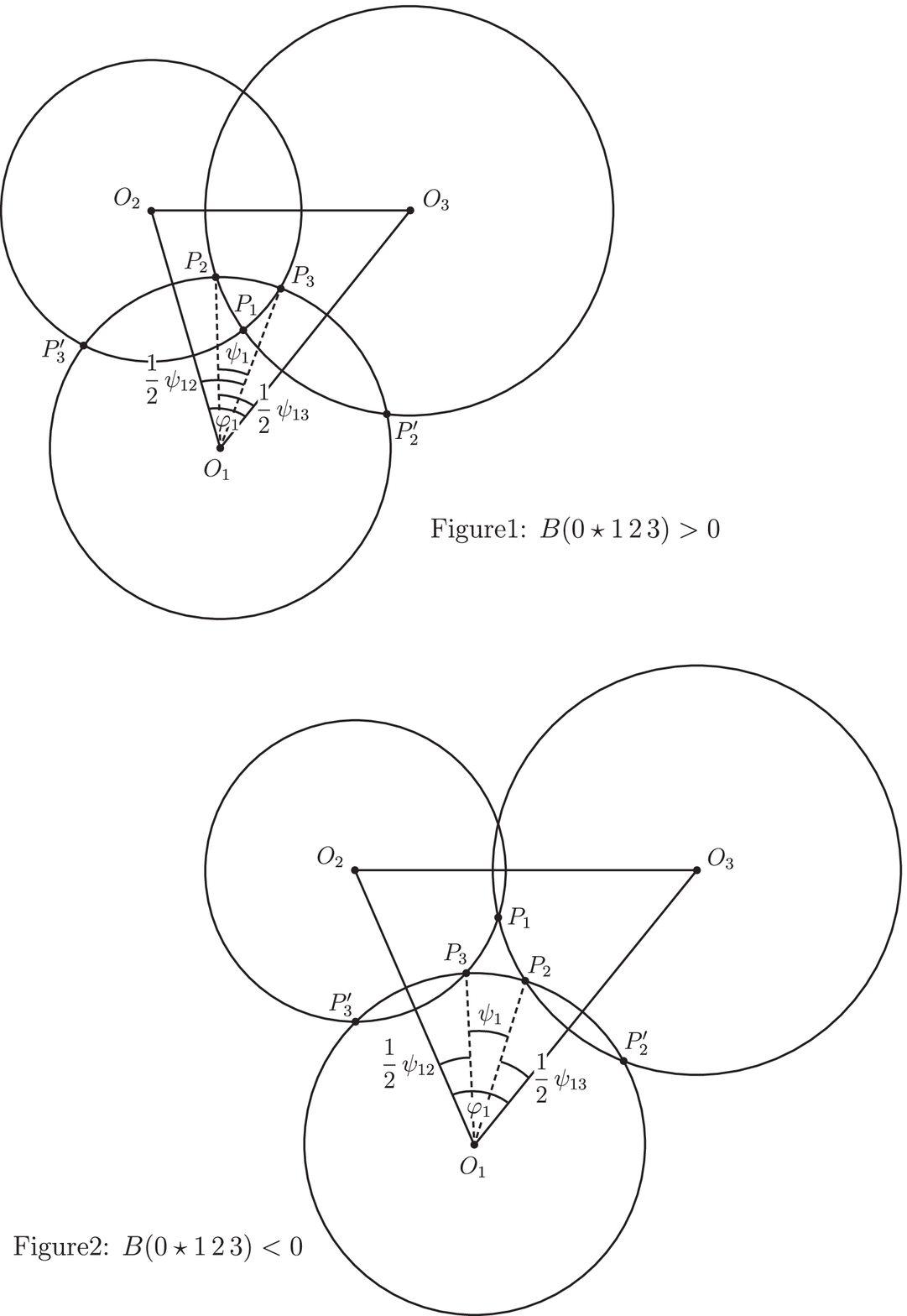}
 \end{center}
\end{figure}
\pagebreak

\bigskip

\section{Restriction to the unit hypersphere}

\bigskip

 We assume further 
  \be
 f_{n+1}(x) = Q(x) - 1,
  \ee
 i.e., $S_{n+1}$ is the unit hypersphere with center $O_{n+1}$
at the origin.

We may assume the linear functions
\be
f_j^\prime(x) :  =  f_j(x) - Q(x) + 1
= \sum_{\nu=1}^n u_{j\nu}x_\nu + u_{j0} \quad(1\le j\le n)
\ee
are normalized such that
the configuration matrix $A^\prime = (a_{jk}^\prime)\ ({0\le j,k\le n})$ of order $n+1$
consisting of 
 \be
 &&a_{j0}^\prime = a_{0j}^\prime =  u_{j0},\\
 &&a_{jk}^\prime = \sum_{\nu=1}^n\, u_{j\nu}u_{k\nu} - u_{j0}u_{k0},
 \ee
satisfies $a_{00}^\prime = -1,\,a_{jj}^\prime = 1\ (1\le j\le n).$
We put further 
\be
f_{n+1}^\prime = 1 - Q(x).
\ee

For the set of indices $J = \{j_1,\ldots, j_p\},\; K = \{k_1,\ldots,k_p\}\subset
 \{0,1,\ldots,n,n+1\}$,
we denote by $A^\prime\!
\left(
\begin{array}{ccc}
 J \\
  K
\end{array}
\right)
$
the subdeterminant with the $j_1,\ldots,j_p$th rows
and the $k_1,\ldots,k_p$th columns. In particular,
we abbreviate $A^\prime\!
\left(
\begin{array}{ccc}
 J  \\
  J  \\
\end{array}
\right)
$
 by $A^\prime(J)$.

The family of 
the hyperplanes $H_j :  f_j^\prime(x) = 0$  define the arrangement of hyperplanes 
${\cal A}^\prime = \bigcup_{j=1}^n H_j$
which correspond to ${\cal A}= \bigcup_{j=1}^n S_j,\ S_j:f_j(x)=0,$ one-to-one. 

The components of the matrix $A^\prime $ are described by the Cayley-Menger determinants 
as follows:
\beq
&&a_{j0}^\prime = \frac{B\!
\left(
\begin{array}{ccc}
 0 &j   &n+1   \\
 0 &\star   &n+1  
\end{array}
\right)
}{\sqrt{-B(0\star j\,n+1)}},\\
&&
a_{jk}^\prime = \frac
{-B\!
\left(
\begin{array}{cccc}
  0&\star   &j& n+1   \\
  0&\star   &k&n+1   \\
\end{array}
\right)
}
{\sqrt{B(0\star j\,n+1)\,B(0\star k\,n+1)}}.
\eeq

$H_j$ has the same intersection
with $S_{n+1}$ as the intersection $S_j\cap S_{n+1}.$

From now on, we shall assume the condition ${ \rm({\cal H}1)}$.

${ \rm({\cal H}1)}$ can be rephrased in terms of the minors of $A^\prime$
as follows:
\be
{ \rm({\cal H}1)}\quad   A^\prime(0\, J) < 0\ (J \subset \partial_{n+1}N),\quad
 A^\prime(J) > 0 \ (1\le |J| , J \subset \partial_{n+1}N).
\ee
Remark that it always holds: $- A^\prime(0\,J) > A^\prime(J) > 0$.

Since $S_{n+1}$ is the unit hypersphere, we have the identity
\be
B(0\star \,n+1) = 2,\, B(0\star j\,n+1) = -1,
\ee
so that
\be
a_{jk}^\prime = - B(0\star j\, k\,n+1) = - \cos\langle j,k\rangle,
\ee
where $\langle j,k\rangle$ denotes the angle subtended by $S_j,S_k$ in $S_{n+1}$.

Let $D = D_{12\ldots\, n+1}^-$ be the (non-empty) real $n$ dimensional domain defined by
\be
D_{12\ldots\,n+1}^- = \bigcap_{j=1}^{n+1}\,D_j^-, \;\;\;
D_j^-  :  f_j^\prime\le0 \subset {\bf R}^n\quad(1\le j\le n+1).
\ee

Then, for any $J \subset \partial_{n+1}N$ such that $|J|=p, \,1\le p\le n - 1,$ 
the intersection $S_{J\,n+1}=S_{n+1} \cap \bigcap_{j\in J} S_j$ 
defines an $n-p-1$ dimensional sphere. In particular,
$\bigcap_{k\in \partial_j \partial_{n+1}N}\,S_k$  consists of two points.

The orientation of ${\bf R}^n$  and $D$ is determined such that
the standard $n$-form $\varpi$ is positive:
\be
\varpi = dx_1\wedge\cdots\wedge dx_{n} > 0.
\ee
We can define the standard volume form on $S_{n+1}$
as
\be
&&\varpi_{n+1} : = \sum_{\nu=1}^n (-1)^\nu \,x_\nu dx_1\wedge\cdots\widehat{dx_\nu}\cdots\wedge dx_n
= 2\,\bigl[\frac{\varpi}{df_{n+1}^\prime}\bigr]
_{S_{n+1}}.
\ee

\bigskip

Let $\Phi^\prime(x)$ be the multiplicative function
\be
\Phi^\prime(x) = \prod_{j\in \partial_{n+1}N} f_j^\prime(x)^{\lambda_j}\quad(\lambda_j \in {\bf R}_{\ge0}).
\ee
We take the value of the many valued function $\Phi^\prime(x)$ such that $\Phi^\prime(x) > 0$ 
at the infinity in ${\bf R}^n.$  

Denote the twisted rational de Rham $(n-1)$-cohomology by
$H_\nabla^{n-1}(X,\Omega^\cdot(*S)) $ and  its dual by $H_{n-1}(X, {\cal L}^*),$
where ${\cal L}^*$ denotes the dual local system on the space
$X =  S_{n+1} - \bigcup_{j \in {\cal A}}S_j$
associated with $\Phi^\prime$.
The covariant differentiation $\nabla $ is given by
\be
\nabla \psi = d\psi + d\log\Phi^\prime\wedge\,\psi.
\ee

The corresponding integral can be expressed as the pairing
\be
H_\nabla^{n-1}(X,\Omega^\cdot(*S))\times H_{n-1}(X, {\cal L}^*) \ni (\varphi  , \mathfrak{z} ) \longrightarrow {\cal J}_\lambda^\prime(\varphi) 
= \int_{\mathfrak{z}} \Phi^\prime(x)\varphi(x)\,\varpi_{n+1}.
\ee
for $\varphi\varpi \in \Omega^{n-1}(*S)$ and  a twisted $(n-1)$-cycle $\mathfrak{z}.$

The following has been proved in [4].

\begin{prop}
$H_\nabla^{n-1}(X,\Omega^\cdot(*S))$ is of dimension $2^{n}$
and has a basis
\be
F_J^\prime = \frac{\varpi_{n+1}}{f_J'},
\ee
where $f_J^\prime$ means the product $\prod_{j\in J}\,f_j^\prime$
and $J$ ranges over the family of all unordered subsets of indices such that $J \subset \partial_{n+1}N$
including the empty set $\emptyset$.
\end{prop}

From now on, we choose a twisted cycle $(n-1)$-cycle $\mathfrak{z}$ such that
\be
\int_{\mathfrak{z}} \Phi^\prime(x)\,\varphi \,\varpi_{n+1}
=\int_{D_{12\ldots\,n+1}^-} |\Phi^\prime(x)|\,\varphi \,\varpi_{n+1}\quad(\varphi \,\varpi_{n+1} \in \Omega^{n-1} (*S)).
\ee
$F_\emptyset^\prime$ means $\varpi_{n+1}$, and we define
\be
{\cal J}_\lambda^\prime(\varphi) = \int_{\mathfrak{z}}\Phi^\prime(x)\,\varphi\,\varpi_{n+1}.
\ee

The derivation of the integral ${\cal J}_\lambda^\prime(\varphi)$ with respect to the parameters
$a_{jk}^\prime, a_{j0}^\prime$ can be expressed  as
\beq
&&d_{A^\prime} {\cal J}_\lambda^\prime(\varphi) = \sum_{j=1}^n da_{j0}^\prime \frac{\partial}{\partial{a_{j0}^\prime}} {\cal J}_\lambda^\prime(\varphi)
+  \sum_{1\le j,k\le n} da_{jk}^\prime
\frac{\partial}{\partial{a_{jk}^\prime}}{\cal J}_\lambda^\prime(\varphi)\nonumber\\
&&{}\qquad \qquad = \int_{\mathfrak{z}} \Phi^\prime(x)\,\nabla_{A^\prime}(\varphi\,\varpi_{n+1}),
\eeq
 where 
  \be
 \nabla_{A^\prime}(\varphi\,\varpi_{n+1}) = d_{A^\prime}(\varphi\,\varpi_{n+1} )+ d_{A^\prime}\log\Phi^\prime(x)\wedge\varphi\,\varpi_{n+1}.
 \ee
 
In addition to the above basis,  it is convenient to introduce the following basis
which we call \lq\lq of second kind":

 \begin{df}
 \be
 F_{*,J}^\prime :  =  F_J^\prime + \sum_{\nu\in J} \frac{A^\prime\!
\left(
\begin{array}{ccc}
 0 &\partial_\nu J     \\
 \nu &\partial_\nu J    
\end{array}
\right)
}{A^\prime(J)}\,F_{\partial_\nu J}^\prime.
 \ee
 In particular, $F_{*,\emptyset}^\prime = F_\emptyset^\prime = \varpi_{n+1}$.
 \end{df}
 
 The differential  $1$-forms defined below will play an essential role in the sequel.
 
 \begin{df}
 \be
 &&\theta^\prime_j  : =  da_{j0}^\prime,\\
&&\theta_{jk}^\prime : =   da_{jk}^\prime - \frac{A^\prime\!
\left(
\begin{array}{ccc}
 0 &k     \\
  j&k      \\  
\end{array}\right)
}{A^\prime(0k)}\,da_{k0}^\prime
- \frac{A^\prime\!
\left(
\begin{array}{ccc}
 0 &j     \\
  k&j      \\  
\end{array}\right)
}{A^\prime(0j)}\,da_{j0}^\prime.
 \ee
 General $\theta_J^\prime$ for $|J| \ge 3$ are defined by induction:
  \be
\theta_J^\prime  : = - \sum_{\nu \in J} \frac{A^\prime\!
\left(
\begin{array}{ccc}
  0& \partial_\nu J   \\
  \nu&\partial_\nu J      \\
\end{array}
\right)
}{A^\prime(0\,\partial_\nu J)} \,\theta_{\partial_\nu J}^\prime \quad(3\le |J|\le n).
 \ee
 \end{df}
 
\medskip

 Denote $\lambda_\infty^\prime = \sum_{j=1}^n \,\lambda_j$
 and $J = \{j_1,\ldots,j_p\},\,|J| = p.$
 
 The following fact has been proved in [4].
 
 \begin{prop}
 The following variation formula holds{\rm :}
  \beq
 \nabla_{A^\prime}(F_\emptyset^\prime) \sim
 \sum_{p=1}^n \sum_{1\le j_1<\cdots< j_p\le n} 
 \frac{\lambda_{j_1}\cdots\lambda_{j_p}}{\prod_{q=1}^{p-1} (-\lambda_\infty - n + q +1)}
 (-1)^p \,\theta^\prime_J \frac{A^\prime(J)}{A^\prime(0\,J)} \,F_{^*,J}^\prime.
 \eeq
 {\rm(}The formula {\rm (4.12)} in {\rm[4]} has  an error. 
  In the {\rm RHS}, the sign $(-1)^p$ should be added as above to the original formula{\rm)}. 
 \end{prop}
 
 For example,
  \be
&&\bullet\;  \nabla_{A^\prime}(F_\emptyset^\prime) \sim - \lambda_1\,\frac{1}{A^\prime(0\,1)}\,\theta_1^\prime\,(F_1^\prime + a_{10}^\prime \,F_\emptyset^\prime)
  - \lambda_2\,\frac{1}{A^\prime(0\,2)}\,\theta_2^\prime\,(F_2^\prime + a_{20}^\prime \,F_\emptyset^\prime)\\
&&{}\qquad \qquad  - \frac{\lambda_1\lambda_2}{\lambda_\infty}\,\frac{A^\prime(1\,2)}{A^\prime(0\,1\,2)}\,\theta_{123}^\prime\,F_{*,\,12}^\prime
\quad (n = 2),\\
&&\bullet\;  \nabla_{A^\prime}(F_\emptyset^\prime) \sim -\sum_{j=1}^3 \lambda_j da_{j0}^\prime F_{*,\,j}^\prime
- \sum_{1\le j<k\le 3}\,\frac{\lambda_j\,\lambda_k}{\lambda_\infty + 1}\frac{A^\prime(jk)}{A^\prime(0jk)}\,\theta_{jk}^\prime \,F_{*\,,jk}^\prime\\ 
&&{}\qquad \qquad - \frac{\lambda_1\lambda_2\lambda_3}{\lambda_\infty(\lambda_\infty+1) }\,\theta_{123}^\prime\,\frac{A^\prime(123)}{A^\prime(0123)} \,F_{*,\,123}^\prime
\quad(n=3).
 \ee

\bigskip

 \section{Analog of Schl\"afli formula (errata and some comments)}
 
\bigskip

The variational formula for the volume of a spherically faced simplex 
in the unit hypersphere was presented  in [4].
In addition to the formulae stated in Theorem I and II here, they make
 a completely integrable system. 

However, some formulae stated there have a few errors.
In this section, we present a correct version as in Theorem III. 
 
Let $P_j\,(1\le j\le n)$ be the points in ${\bf R}^n$ such that
\be
\{P_j\}= \bigcap_{k\in \partial_j N} H_k\cap S_{n+1}.
\ee

We can take the Euclidean coordinates $x_1,\ldots,x_n$
such that the polynomials $f_j$ have the following expressions:
\beq
f_j^\prime(x) = \sum_{\nu=1}^{n+1-j} \,u_{j\,\nu} x_\nu + u_{j0}\quad(1\le j\le n).
\eeq

We assume for simplicity that $u_{j\,n+1-j} = 2\alpha_{j\,n+1-j}> 0\ (1\le j\le n)$
and that $P_j$ satisfies the condition (${\cal H}$2).

We have the equalities
\beq
\prod_{j=p+1}^{n} u_{j\,n- j +1} =  \sqrt{-\,A^\prime(0\,p+1\,\ldots n)}\quad(1\le p\le n ).
\eeq

The affine subspace $\bigcap_{j= n-p+1}^n \,H_j$ contains the $n-p-1$ dimensional sphere
$S_{n-p+1\ldots n\,n+1} = \bigcap_{j= n-p+1}^n S_j \cap S_{n+1}$ with radius
\be
r_{n-p+1\ldots n\,n+1} = \sqrt{-\frac{A^\prime(n-p+1\ldots n)}{A^\prime(0\,n-p+1\ldots\,n) }}.
\ee
Denote by $\widetilde{\Delta}[P_1,P_2,\ldots,P_{n}]$
be the pseudo $n-1$-simplex in $S_{n+1}$ with hyperspherical 
faces with vertices $P_j$ such that their sign of orientation is $(-1)^{\frac{n(n-1)}{2}}.$ 
The support of $\widetilde{\Delta}[P_1,P_2,\ldots,P_{n}]$ coincides with $D = D_{12\ldots\, n+1}^-.$

By definition, the following properties are valid.

\begin{lm}

{\rm (i)} 
\be
&&df_n^\prime\wedge\cdots \wedge df_1^\prime>0
\ee
on $D $.

{\rm(ii)}
The pseudo $(n-1)$-simplex $\widetilde{\Delta}[ P_1,P_2,\ldots,P_{n}]$ has the sign    
$(-1)^{\frac{n(n-1)}{2}}$ of orientation such that 
\be
\widetilde{\Delta}[ P_1,P_2,\ldots,P_{n}] = (-1)^{\frac{n(n-1)}{2}} \, S_{n+1}\cap D.
\ee
\end{lm}

\begin{pf}
Indeed, we can show that
\beq
df_{n}^\prime\wedge\cdots\wedge df_1^\prime = \prod_{j=1}^n \,u_{j\,n-j+1}\,\varpi > 0.
\eeq

(ii) follows from $({\cal H}2).$
\qed
\end{pf}

\medskip

Let $v_\emptyset$ be the volume of the pseudo $(n-1)$-simplex $\widetilde{\Delta}[ P_1,P_2,\ldots,P_{n}]$
defined by
\be
v_\emptyset = \int_{\widetilde{\Delta}[ P_1,P_2,\ldots,P_{n}]} \varpi_{n+1}> 0,
\ee 
where the orientation of $\widetilde{\Delta}[ P_1,P_2,\ldots,P_{n}]$ is chosen such that
$\varpi_{n+1}$ should be positive on it.

We are interested in the variation 
formula for $v_{\emptyset}$, which can be expressed in terms of the lower dimensional
volumes of the faces of $\widetilde{\Delta}[ P_1,P_2,\ldots,P_{n}]$.

Every face of the pseudo simplex is included in some $S_{J\,n+1}.$
The  $S_{J\,n+1}$ is defined as an $n-p-1$ dimensional sphere
with radius 
\be
r_{J\, n+1} = \sqrt{- \frac{A^\prime(J)}{A^\prime(0\,J)}}.
\ee

We can consider the $n-p-1$ dimensional volume $v_J\,(|J| = p)$ relative to  the corresponding standard
volume form $\varpi_{J\,n+1}^\prime$ on the $n-p-1$ dimensional sphere:
\be
&&v_J = \int_{\widetilde{\Delta}[ P_1,P_2,\ldots,P_{n}] \cap S_J}\, |\varpi_{J\,n+1}^\prime|,\\
\ee
where
\be
|\varpi_{J\,n+1}^\prime|
= r_{J\,n+1}^{n-p-1}\,|\varpi_{J\,n+1}| > 0.
\ee
The orientation of $\widetilde{\Delta}[ P_1,P_2,\ldots,P_{n}] \cap S_J$ is chosen such that
$\varpi_{J\,n+1}^\prime$ should be positive :  $\varpi_{J\,n+1}^\prime = |\varpi_{J\,n+1}^\prime|$,
\ $|\varpi_{J\,n+1}^\prime|$ being the absolute value of $\varpi_{J\,n+1}^\prime$.

When $J = \{n-p+1\ldots n\},$ we can give an explicit expression
for $\varpi_{n-p+1\ldots n\,n+1}$ as follows:
\be
&&f_j^\prime(x) = 0 \quad(n-p+1\le j\le n+1),\\
&& \sum_{j= p+1}^n \,x_j^2 = r_{n-p+1\ldots n+1}^2,
\ee
where 
\be
r_{n-p+1\ldots n+1} = \sqrt{ - \frac{A^\prime(n-p+1\ldots n)}{A^\prime(0\,n-p+1\ldots n)} } .
\ee

The standard volume form on $S_{n-p+1\,\ldots\,n+1}$ is given by
\beq
&&\varpi_{n-p+1\ldots n+1}^\prime 
 = \sum_{\nu= p+1}^n \,(-1)^\nu \,
\frac{
x_\nu\,dx_{p+1}\wedge\cdots\widehat{dx_\nu}\cdots\wedge dx_n
}
{r_{n-p+1\ldots n+1}}\nonumber\\
&&{}\qquad \qquad \quad = r_{n-p+1\ldots n+1}^{n-p-1}\,\varpi_{n-p+1\ldots n+1},
\eeq
where
\be
\varpi_{n-p+1\ldots n+1} = \sum_{\nu= p+1}^n \,(-1)^\nu \,
\xi_\nu\,d\xi_{p+1}\wedge\cdots\widehat{d\xi_\nu}\cdots\wedge d\xi_n.
\ee
through the transformation
\be
&&x_\nu = r_{n-p+1\ldots n+1}\,\xi_\nu\quad(p+1\le \nu\le n),
\ee
such that $\sum_{\nu = p+1}^n \,\xi_\nu^2 = 1.$

The following Lemma follows by definition of residue formula.

\begin{lm}
For $J = \{ j_1,\ldots,j_p\}\,(1\le  j_1<\ldots < j_p\le n),$
\be
\bigl[\frac{\varpi_{n+1}}
{df_{j_p}^\prime\wedge\cdots\wedge df_{j_1}^\prime}\bigr]_{S_{j_1\,\ldots \,j_p}}
= \frac{1}{\sqrt{A^\prime(J)}}\,
 \varpi_{J\,n+1}^\prime.
\ee
In particlular,
\be
\bigl[\frac{\varpi_{n+1}}
{df_n^\prime\wedge\cdots\widehat{df_j^\prime}\cdots\wedge df_{1}^\prime}\bigr]_{P_j}
=  \frac{(-1)^{n-j}}{\sqrt{A^\prime(\partial_j \partial_{n+1} N})}\quad(1\le j\le n),
\ee
since  $[f_j^\prime]_{P_j}$ at the point $P_j$ of $S_{\partial_j\,\partial_{n+1} N}\cap D$
is negative. 
\end{lm}

\begin{pf}
To prove Lemma 16, we may assume that 
$j_1=n-p+1,\ldots, j_p=n$ 
and $f_j'$ are represented by the reduced form (42). 
A direct calculation and (43) show the following identity 
\be
&&d(1-Q(x))\wedge df_n'\wedge \cdots \wedge df_{n-p+1}'\wedge 
\sum_{\nu=p+1}^n (-1)^{\nu}x_{\nu}dx_{p+1}\wedge \cdots 
\widehat{dx_{\nu}}\cdots \wedge dx_n \\
&&{}\qquad \qquad =2\prod_{q=1}^p u_{n-q+1\ q}\ (\sum_{\nu=p+1}^n x_{\nu}^2)\
 \varpi \\
&&{}\qquad \qquad =2\sqrt{-A'(0\ n-p+1\cdots n)}\ r_{n-p+1\ldots n+1}^2\
 \varpi.
\ee
Hence,
\be
&&[df_n'\wedge \cdots \wedge df_{n-p+1}'\wedge 
\sum_{\nu=p+1}^n (-1)^{\nu}x_{\nu}dx_{p+1}\wedge \cdots 
\widehat{dx_{\nu}}\cdots \wedge dx_n]_{S_{n+1}} \\
&&{}\qquad \quad =\sqrt{-A'(0\ n-p+1\cdots n)}\ r_{n-p+1\ldots n+1}^2
\varpi_{n+1}.
\ee
Namely,
\be
&&\varpi_{n-p+1\ldots n+1}^{\prime}=\frac{\sum_{\nu=p+1}^n (-1)^{\nu}x_{\nu}dx_{p+1}\wedge \cdots 
\widehat{dx_{\nu}}\cdots \wedge dx_n}{r_{n-p+1\ldots n+1}} \\
&&{}\qquad =\sqrt{A'(n-p+1\cdots n)} 
\ [\frac{\varpi_{n+1}}{df_n'\wedge \cdots \wedge df_{n-p+1}'}]_{S_{n-p+1\ldots n+1}}.
\ee

General volume forms $\varpi_{J\ n+1}^{\prime}$ can be explicitly written 
by the use of suitable coordinates transformed by isometry.
\qed
\end{pf}

\medskip

The next Theorem has been essentially stated in [4] (Theorem 8), 
but has some errors in the formulae (5.6) therein.
Here we state a correct version,
which follows from Proposition 14.

\bigskip

\bigskip


[{\bf Theorem III}]

\medskip

{\it
For $v_\emptyset = v(D_{12\ldots\,n}^-)$, we have 
\beq
&&d_{A^\prime}v_\emptyset =  - \sum_{p=1}^{n-1} \sum_{|J| = p}\,(-1)^p\,\frac{(n-p-1)!}{(n-2)!}\,\theta_J^\prime
\,\frac{\sqrt{A^\prime(J)}}{A^\prime(0\,J)} \,v_J\nonumber\\
&&{}\qquad \quad + (-1)^n \,\frac{1}{(n-2)!}\,\frac{1}{\sqrt{- A^\prime(0\,1,\ldots, n)}}\,\theta_{12\ldots\,n}^\prime,
\eeq
where $J$ ranges over the collection of unordered subsets of $\{1,2,\ldots, n\}$ 
and $|J| = p.$
}

\bigskip

To prove this Theorem, we need the following Lemma equivalent to Proposition 6.

\begin{lm}
We have the identity
\be
&&\bigl[\frac{1}{f_j^\prime}\bigr]_{P_j} = \bigl[\frac{1}{f_j}\bigr]_{P_j}\\
&&{}\quad  = \frac{
\sqrt{- A^\prime(\partial_j\partial_{n+1}N)\,A^\prime(0\,\partial_{n+1}N)} 
+ A^\prime\!
\left(
\begin{array}{ccc}
  0&\partial_j\partial_{n+1}N    \\
  j&  \partial_j\partial_{n+1}N     \\
\end{array}
\right)
}{- A^\prime(\partial_{n+1}N)} < 0,
\ee
so that
\be
\bigl[\frac{1}{f_j^\prime}\bigr]_{P_j} + \frac{A^\prime\!
\left(
\begin{array}{ccc}
  0&\partial_j\partial_{n+1}N    \\
  j&  \partial_j\partial_{n+1}N     \\
\end{array}
\right)}{A^\prime(\partial_{n+1}N)}
=  - \frac{\sqrt{- A^\prime(\partial_j\partial_{n+1}N)\,A^\prime(0\,\partial_{n+1}N)}}{A^\prime(\partial_{n+1}N)}.
\ee
\end{lm}


\bigskip

{\bf Proof of Theorem III.}

\medskip

Take $\lambda_j$ such that all $\lambda_j = \varepsilon > 0$ in the formula
{\rm(41)}. Then 
{\rm(40)} shows that
\be
&&d_{A^\prime}\,V_\emptyset = \lim_{\varepsilon\downarrow0}\,d_{A^\prime}\,\int_{\widetilde{\Delta}[P_1,\ldots, P_n]}|\Phi^\prime(x)|\,\varpi_{n+1}\\
&&{}\qquad \quad = \lim_{\varepsilon\downarrow0}\int_{\mathfrak{z}}\Phi^\prime(x)\,\nabla_{A^\prime}(\varpi_{n+1})\\
&&{}\qquad \quad =\lim_{\varepsilon\downarrow0}\int_{\widetilde{\Delta}[P_1,\ldots, P_n]} |\Phi^\prime(x)|\, \nabla_{A^\prime} (\varpi_{n+1}).
\ee

In view of the formula {\rm(4.11)} and the proof of {\rm Theorem 7} in {\rm[3]}, we have 
only to check the following fact:
\beq
&&\lim_{\varepsilon\downarrow0} \, 
\frac{\prod_{j=1}^n \lambda_j}
{\prod_{q=1}^{n-1} (-\lambda_\infty - n+q+1)}
\frac{A^\prime(\partial_{n+1}N)}{A^\prime(0\,\partial_{n+1}N)}
\Bigl\{{\cal J}_\lambda(\frac{1}{f_{\partial_{n+1}N }})\Bigr.
\nonumber\\
&&\Bigl.
 + \sum_{j\in \partial_{n+1}N}
\frac{A^\prime\!
\left(
\begin{array}{ccc}
 0 &\partial_j\partial_{n+1}N    \\
 j &   \partial_j\partial_{n+1} N  
\end{array}
\right)
}{A^\prime(\partial_j\partial_{n+1}N)}
{\cal J}_\lambda\,
(\frac{1}{
f_{\partial_j\partial_{n+1}N}
})
\Bigr\}
 = \frac{(-1)^n}{(n-2)! \sqrt{- A^\prime(01\ldots n)}}.\nonumber
 \\
\eeq
 
By the residue theorem, the {\rm LHS} reduces to $n$ pieces of  point
measures at $P_j$ and equals
\be
&&\lim_{\varepsilon\downarrow0} \frac{\varepsilon^n}{\prod_{q=1}^{n-1} (-n\varepsilon - n + q +1)}
  \,{\cal J}(F_{^*,\partial_{n+1}N}^\prime)\,\frac{A^\prime(\partial_{n+1}N)}{A^\prime(0\,\partial_{n+1}N)}\,
\\
&&= - \lim_{\varepsilon\downarrow0} \frac{\varepsilon^{n-1}}{n\,\prod_{q=1}^{n-2} (-n\varepsilon - n + q +1)}
  \,{\cal J}(F_{^*,\partial_{n+1}N}^\prime)\,\frac{A^\prime(\partial_{n+1}N)}{A^\prime(0\,\partial_{n+1}N)}\,
\\
&&= \sum_{j=1}^n \frac{(-1)^{n-1}}{n(n-2)!}
\{\bigl[\frac{1}{f_j^\prime}\bigr]_{P_j} + 
  \frac{
A^\prime\!
\left(
\begin{array}{ccc}
0&\partial_j\partial_{n+1}N\\
  j&\partial_j\partial_{n+1}N     \\
\end{array}
\right)
}{A^\prime(\partial_{n+1}N)}\}\,\frac{A^\prime(\partial_{n+1}N)}{A^\prime(0\,\partial_{n+1}N)}\,
\,
\Bigl|\bigl[\frac{\varpi_{n+1}}{df_n^\prime\wedge\cdots\widehat{df_j^\prime}\cdots\wedge df_1^\prime}\bigr]_{P_j}\Bigr|.
\ee

 On the other hand, we have 
 \be
 \bigl[\frac{\varpi_{n+1}}{df_n^\prime\wedge\cdots\widehat{df_j^\prime}\cdots\wedge df_1^\prime}\bigr]_{P_j}
 = \frac{(-1)^{n+1-j}}{\sqrt{A^\prime(\partial_j\partial_{n+1}N)}}.
 \ee
 
 Each term in the summand of the {\rm RHS} does not depend on $j$
and is equal to 
\be
 \frac{(- 1)^{n-1}}{n\,(n-2)!}
\{ \frac{
\sqrt{- A^\prime(\partial_j\partial_{n+1}N)\,A^\prime(0\,\partial_{n+1}N)} 
}{A^\prime(\partial_{n+1}N)}
\}\,\frac{1}{\sqrt{A^\prime(\partial_j\partial_{n+1}\,N)}}
=
 \frac{(- 1)^{n-1}}{n\,(n-2)!}\,
\frac{
\sqrt{- A^\prime(0\,\partial_{n+1}N)} 
}{A^\prime(\partial_{n+1}N)}.
\ee
 Hence, the {\rm LHS} of {\rm(47)} becomes
 
 \be
&&\lim_{\varepsilon\downarrow0} \frac{\varepsilon^n}{\prod_{q=1}^{n-1} (-n\varepsilon - n + q +1)}
 \frac{A^\prime(\partial_{n+1}N)}{A^\prime(0\,\partial_{n+1}N)} \,{\cal J}(F_{^*,\partial_{n+1}N}^\prime)\\
&&{}\qquad \qquad  =
\frac{(-1)^n}{(n-2)!}\,\frac{1}{\sqrt{-A^\prime(0\,\partial_{n+1}N})}\,\theta_{\partial_{n+1}N}^\prime.
 \ee
 In this way, we have proved {\rm Theorem III}.
 \qed
 
\medskip

\begin{re}
In three dimensional case, i.e., for $n = 3$, 
$D_{123}^- $ is a pseudo triangle $\tilde{\Delta}P_1P_2P_3$ with circular arc sides. Theorem III shows the identity
\beq
&&d_{A^\prime} \,v_\emptyset = \sum_{j=1}^3 \theta_{j}^\prime \,\frac{1}{A^\prime(0\,j)}\, v_j 
- \sum_{j<k} \, \theta_{jk}^\prime\, \frac{\sqrt{A^\prime(jk)}}{A^\prime(0\,j\,k)}\nonumber\\
&&{}\quad \quad  - \frac{1}{\sqrt{- A^\prime(0\,1\,2\,3)}} \,\theta_{123}^\prime.
\eeq

On the other hand, Gauss-Bonnet theorem shows the identity
\beq
v_\emptyset = 2\pi - \sum_{j=1}^3 \,a_{j0}^\prime\,v_j - \sum_{j<k} (\pi - \langle jk \rangle), 
\eeq
where $\langle jk\rangle$\, denotes the angle of the triangle at $P_l\,(\{j,k,l\} 
:\,{\rm a\,permutation\,of}\\
 \, \{1,2,3\})$ such that
\be
a_{jk}^\prime = - \cos\langle jk\rangle,
\ee
and $a_{j0}^\prime$ is the geodesic curvature of the arc $\partial D_{123}^-\cap S_j.$

We can see by a direct calculation that the differential of (49) coincides with (48).
 Gauss-Bonnet theorem was extended into a higher dimensional polyhedral domain by Allendoerfer-Weil 
 (see the second formula in [2]).
However, in the case of a spherically faced simplex,
the formula (46) does not seem to generally coincide with  the differential  of the latter.  

\end{re}


\bigskip

\bigskip

\noindent
{\Large{\bf Appendix}}\quad {\Large{\bf Elementary proof of Theorem II (i)}}

\bigskip

Denote by $P_j\,(1\le j\le n+1))$ the vertex points of the $n$-simplex $D_N^+$
such that $P_j \in  \partial D_N^+\cap \bigcap_{k\ne j} S_k$.
For the ordered set $J = \{j_1,\ldots, j_p)\subset N$
such that $j_1 > j_2 > \cdots > j_p,\ |J| = p$  
and $J^c  = \{j_1^*> \cdots > j_{p}^* \},$
$\tilde{\Delta}O_J\,P_ {J^c}$ means the $n$-cell
\be
\tilde{\Delta}O_{J}\,P_ {J^c} = \tilde{\Delta}O_{j_1}\ldots O_{ j_{p}}\,P_{j_{1}^*}\,\ldots \,P_{j_{n-p+1}^*}
\ee
with the vertices $O_{j_1},\ldots, O_{j_{p}}$ and 
$P_{j_1^*},\ldots, P_{j_{n-p+1}^*}$.
Notice that $\tilde{\Delta}P_{j_1}^*\,\ldots\,P_{j_{n-p+1}^*}$ $= S_J\cap D_N^+$ are $(p-1)$-pseudo simplex
with the faces  $S_k\cap S_J\cap D_N^+\, (k\in J^c)$
in the $n-p$
dimensional sphere $S_J = \bigcap_{j\in J} S_j .$

We have the cell decomposition of ${\Delta}O_{n+1}\,\ldots \,O_2O_1:$
\be
{\Delta}O_{n+1}\,\ldots \,O_2O_1 = - \sum_{p=0}^n\, \sum_{|J| = p} \,\tilde{\Delta}O_{J}\,P_{J^c}\,\varepsilon_J,
\ee
where $\varepsilon_J$ denotes $(-1)^{\sum_{j\in J^c} j}\,(-1)^{\frac{(n-p)(n-p+1)}{2}}$.

Hence we have the identity for their volumes 
\be
v(\Delta \,O_{n+1}\ldots \,O_1) = 
\sum_{J \subset N ,\,|J|\le n} v(\Delta O_J P_{J^c}),
\ee
or equivalently,
\be
v(\tilde{\Delta}P_N) = v(\Delta O_{n+1}\,\ldots\,O_1) - \sum_{J \subset N ,\,1\le |J|\le n} v(\Delta O_J P_{J^c}).
\ee

The identity stated in Theorem II (i) is a direct consequence of  the following Lemma.

\medskip

\begin{lm}
\be
v(\tilde{\Delta}O_J\,P_{J^c}) = \frac{(n-p)!}{(n-1)!}\, 
\sqrt{\frac{(-1)^{p+1}\,B(0\star J)}{2^p}}\, v_J.
\ee
\end{lm}

 \begin{pf}
 Without losing generality, we may assume that 
 $f_j$ have the standard form (1), (2) and $J = \{n+1,n,\ldots, n-p+2\}.$

 $O_j \,(n-p+2\le j\le n+1)$ can be expressed as
 \be
 O_j = ( -\alpha_{j1},\,\ldots,\, -\alpha_{j,n-j+1},\,0,\ldots, 0)\quad(\alpha_{j,n-j+1} >  0). 
 \ee

 The spherical $(n-p)$-simplex $\tilde{\Delta}P_{n-p+1}\ldots P_1$
 with support $D_N^+\cap S_J$
 is defined by the equations for $\xi = (\xi_1,\ldots,\xi_n):$
  \beq
 f_j(\xi) = 0\ (n-p+2\le j\le n+1),\ f_k(\xi) \ge 0\ (1\le k\le n-p+1).
 \eeq
The coordinates $\xi_j \ (1\le j\le p-1)$ are uniquely determined by (50)
and denoted by $\gamma_j$.

$\xi $ ranges over the $n-p$ dimensional sphere
\be
 \xi = (\gamma_1, \ldots, \gamma_{p-1}, \xi_p,\ldots,\xi_n)
\ee
under the condition
\be
&& f_k(\xi) \ge 0 \quad  (1\le k\le n-p+1),\\
&&\sum_{j= p}^{n+1}\,\xi_j^2 = r_{n-p+2\,\ldots\,n+1}^2,
\ee
where $r_{n+1\ldots n-p+2}$ denotes the radius of the hypersphere
$S_J$:
\be
r_{ n-p+2\,\ldots \,n+1} = \sqrt{\frac{(-1)^p \,B(0\,n-p+2\,\ldots\,n+1) }{2^{p-1}}}.
\ee

 The $n$-pseudo simplex $\tilde{\Delta}O_{n+1}\ldots O_{n-p+2}P_{n-p+1}\ldots\, P_{1}$
 consists of the union of the $p$-simplex $\Delta O_{n+1}\ldots O_{n-p+2}\xi$
 with $\xi \in \tilde{\Delta}(P_{n-p+1}\ldots P_1)$: 
 \be
\tilde{\Delta} O_{n+1}\ldots O_{n-p+2}P_{n-p+1}\ldots\, P_{1}
= \bigcup_{\xi \in \tilde{\Delta}P_{n-p+1}\ldots\, P_{1}}\Delta O_{n+1}\ldots O_{n-p+2}\xi .
 \ee
 
 Namely, every point of $\tilde{\Delta} O_{n+1}\ldots O_{n-p+2}P_{n-p+1}\ldots\, P_{1}$
 is parametrized by the expressions: 
 \be
&& x_j = - \sum_{k=1}^{j} y_k \, \alpha_{n-k+1,j} + y_0 \,\gamma_j\quad(1\le j\le p-1),\\
&& x_j = y_0\,\xi_j\quad(p\le j\le n)
 \ee
 such that $y = (y_0,\ldots, y_{p-1})$ ranges over the $p$-convex set
 \be
\delta_p :   y_j \ge 0\,(0\le j\le p-1),\quad \sum_{j=0}^{p-1}\,y_j \le 1.
 \ee
 
 As a result in view of (3),(9) and (45),
 \be
 && v(\tilde{\Delta} O_{n+1}\ldots O_{n-p+2}P_{n-p+1}\ldots\, P_{1}) \\
&&{}\quad = \int_{\tilde{\Delta} O_{n+1}\ldots O_{n-p+2}P_{n-p+1}\ldots\, P_{1}
} \, |dx_1\wedge\cdots\wedge dx_{p-1}\wedge dx_p\wedge\cdots dx_n|\\
&&{}\quad = \prod_{j=1}^{p-1}\alpha_{n-j+1,j}\,\int_\delta  y_0^{n-p} dy_1\wedge \cdots\wedge dy_{p-1}\wedge dy_0, \\
&&\int_{\tilde{\Delta}P_{n-p+1}\ldots P_{1}}\,|\sum_{\nu=p}^n (-1)^{\nu-p}\,\xi_\nu\,d\xi_p\wedge\cdots<d\xi_\nu>\cdots\wedge d\xi_n|\\
&&
{}\quad = \frac{(n-p)!}{(n-1)!}\,\sqrt{\frac{(-1)^{p+1} \,B(0\star\,n-p+2\ldots n+1)}{2^p}}\, 
v_{n-p+2\ldots n+1}.
 \ee
 In this way, Lemma 18 has been proved in case where $J = \{n+1\,\ldots\, n-p+2\}$. 
 Therefore it also holds true for general $J$ because of symmetry.
 \qed
  \end{pf}
 
 Theorem I (i) can also be proved in a similar way.


\bigskip

{\bf Acknowledgement}\quad
The authors appreciate a useful suggestion due to M.Ito for drawing the figures.
 

\bigskip

\bigskip

\bigskip

\bigskip

\bigskip

\bigskip

\noindent Kazuhiko AOMOTO,\\
5-1307 Hara, Tenpaku-ku, Nagoya-shi, 468-0015, Japan.\\
e-mail: kazuhiko@aba.ne.jp

\bigskip

\noindent Yoshinori MACHIDA,\\
4-9-37 Tsuji, Shimizu-shi, Shizuoka-shi, 424-0806, Japan. \\
e-mail: yomachi212@gmail.com

\end{document}